\documentclass{./ftpfiles/cjour/cjour}
\usepackage{a4}
\usepackage[dvips]{graphicx}
\usepackage{amsmath}
\usepackage{amsfonts}
\usepackage{amssymb}
\usepackage{psfrag}
\usepackage{lscape}
\newcommand{\DS}{\displaystyle}
\newcommand{\TST}{\textstyle}
\newcommand{\SSS}{\scriptstyle}
\newcommand{\el}{\\[2ex]}
\newcommand{\diffx}[2]{\text{$\DS\frac{\partial #1}{\partial #2}$}}
\newcommand{\Diffx}[2]{\text{$\DS\frac{\hbox{d} #1}{\hbox{d} #2}$}}
\newcommand{\dd}[2]{{\partial #1}/{\partial #2}}
\renewcommand{\thetable}{\thesection.\arabic{table}}
\def\dfrac#1#2{\DS\frac{#1}{#2}}
\def\u#1{\mbox{\boldmath $#1$}}
\def\half{\frac12}
\def\eqref#1{(\ref{#1})}
\def\refer#1#2#3%
{\noindent{\hangindent=1cm \hangafter=1{#1\ #2\ #3\hfill}
\par\filbreak}}

\def\point{{\quad .}}
\def\ljump{[\![}
\def\rjump{]\!]}

\def\text#1{\hbox{#1}}
\def\d{\hbox{d}}
\def\V#1{{\mathbf #1}}

\def\dt{\triangle t}
\def\dx{\triangle x}
\def\halb{{\SSS \frac{1}{2}}}

\def\meanW{{\overline{\V w}}}
\def\xft{x_{Ft}}
\def\xtl{x_{Tl}}
\def\sinZ{\sin\!\zeta}
\def\cosZ{\cos\!\zeta}
\def\tanD{\tan\!\delta}
\begin{document}

%%%%% To be entered at Academic Press: =====>>

% \journame{}
% \articlenumber{}
% \yearofpublication{}
% \volume{}
% \cccline{}
% \received{}
% \revised{}
% \accepted{}

\authorrunninghead{Y.C. Tai, S. Noelle, J.M.N.T. Gray
  \& K. Hutter}
\titlerunninghead{Shock-Capturing and Front-Tracking for Granular Avalanches}

% communication line, use: \commline{Communicated by...}
% \commline{ }

%\setcounter{page}{261} %% This command is optional. 

%% <<== End of commands to be entered at Academic Press

%%  Authors, start here ==>>

%\draft % Optional, will cause a line at the bottom of each page
%% with the words `Draft' and the time and date that the article
%% was LaTeXed. Will also double space text.

\title{Shock-Capturing and Front-Tracking Methods
  for Granular Avalanches}

%\subtitle{}

\author{Y.C. Tai${}^1$, S. Noelle${}^2$, J.M.N.T. Gray${}^3$
  \& K. Hutter${}^1$}
\affil{${}^1$Institut f\"ur Mechanik,
    Technische Universit\"at Darmstadt, 64289 Darmstadt, Germany. \\
    tai@mechanik.tu-darmstadt.de, hutter@mechanik.tu-darmstadt.de}
\affil{${}^2$Institut f\"ur Angewandte Mathematik, 
    Universit\"at Bonn, 53115 Bonn, Germany. noelle@iam.uni-bonn.de}
\affil{${}^3$Department of Mathematics,
    University of Manchester, Manchester M13 9PL, UK. ngray@ma.man.ac.uk}

%%%%%%%%%%%%
%% More than one author with separate affiliations, either:
%
%\author{First Author Name$^\dagger$ and Second Author Name$^\ddagger$}
%\affil{$^{\dagger}$First Author Affiliation, $^{\ddagger}Second Author
%       affiliation}

% or

% \author{Author name}
% \affil{Affiliation}
% \and
% \author{Author name}
% \affil{Affiliation}
%%%%%%%%%%%%

%% \thanks command:
%% Can use \thanks{} in title to have footnote number appear and
%% footnote at the bottom of the page. i.e.,
%% \title{This is the title\thanks{Supported by grant no....}}

%% In \authors or \affil, can use \thanks{} to have asterisk, 
%%   dagger or double dagger appear
%%   and text appear at the bottom of the title page. i.e.,

%\authors{D. Adalsteinsson and J. A. Sethian\thanks{Supported in part by the
%Applied Mathematics Subprogram of the...}}

%%%%%%%%%%%%

%\email{address}

%optional
%\dedication{Dedicated to...}
  
  \abstract{Shock formations are observed in granular avalanches when
    supercritical flow merges into a region of subcritical flow.  In
    this paper we employ a shock-capturing numerical scheme for the
    one-dimensional Savage-Hutter theory of granular flow to describe
    this phenomenon. A Lagrangian moving mesh scheme applied to the
    non-conservative form of the equations reproduces smooth solutions
    of these free boundary problems very well, but fails when shocks
    are formed. A non-oscillatory central difference (NOC) scheme with
    TVD limiter or WENO cell reconstruction for the conservative
    equations is therefore introduced. For the avalanche free boundary
    problems it must be combined with a front-tracking method,
    developed here, to properly describe the margin evolution. It is
    found that this NOC scheme combined with the front-tracking module
    reproduces both the shock wave and the smooth solution accurately.
    A piecewise quadratic WENO reconstruction improves the
    smoothness of the solution near local extrema.
    The schemes are checked against exact solutions for
    (1) an upward moving shock wave,
    (2) the motion of a parabolic cap down an inclined plane and
    (3) the motion of a parabolic cap down a curved slope
    ending in a flat run-out region, where a shock is formed
    as the avalanche comes to a halt.}

% text should be lower case, unless caps are necessary for meaning


\begin{article}

% \contents is optional, will make a list of all section heads 
% that appear in the article
%\contents

%optional
% for those that like to start with section zero:
%\zerosection{Introduction}
%
% Section 1
%

\section{Introduction}

Snow avalanches, landslides, rock falls and debris flows are extremely
dangerous and destructive natural phenomena of which the occurance
has increased during the past few decades. Their human impact has
become so significant that the United Nations declared 1990--2000
{\it International Decade for Natural Disasters Reduction (IDNDR)}.
Research on the protection of habitants from floods, debris flows and
avalanches is under way worldwide, and many institutions focus
on the numerical prediction of such flows under ideal as well as
realistic conditions.

One of the models that has become popular in recent years is the
Savage-Hutter (SH) avalanche theory for granular materials
\cite{SavageHutter:1989,SavageHutter:1991}.
In the past decade numerical techniques were developed to solve
the SH-governing differential equations for typical moving boundary
value problems
\cite{SavageHutter:1989,SavageHutter:1991,HutterKoch:1991,
GreveHutter:1993,GreveKochHutter:1994,KochGreveHutter:1994,
HutterSiegelSavageNohguchi:1993,HutterKochPlussSavage:1995,
GrayWielandHutter:1999,WielandGrayHutter:1999,GrayTai:1998A,
GrayTai:1998B}. These
techniques are based on a Lagrangian moving mesh finite-difference scheme in which
the granular material is divided into quadrilateral cells (2D) or
triangular prisms with flat tops (3D). Exact similarity
solutions of the SH-equations were constructed in spatially
one-dimensional chute flows
\cite{SavageHutter:1989,SavageNohguchi:1988,NohguchiHutterSavage:1989}
and for two-dimensional unconfined flows
\cite{HutterNohguchi:1990,HutterGreve:1993}.
In the case of chute flows it was shown that the solutions obtained by the Lagrangian
integration procedure approximate the exact parabolic similarity solution
very accurately, and these theoretical and numerical results are in good agreement with
experimental avalanche data. Similar agreement between theoretical, numerical and
experimental data was also obtained for the two-dimensional flow configurations
(cf. above references). In these Lagrangian schemes explicit artificial numerical
diffusion was incorporated to maintain stability. In doing so the
quality of resolution deteriorates. In fact, the adequacy of these
numerical solutions can be challenged because of uncontrolled
spreading due to this diffusion. It was also observed that the
Lagrangian schemes loose their stability (or else unjustified
artificial diffusion must be applied) whenever internal shocks are
formed. This appears to occur whenever the avalanche moves from an
extending to a contracting flow configuration. These shocks are
travelling waves which form bumps with steep gradients on the free
surface, which is thicker on the down slope side.
It is therefore natural to develop conservative high-resolution
shock-capturing numerical techniques that are able to resolve the
steep surface gradients and identify the shocks often observed in
experiments but not captured by the Lagrangian finite difference
scheme.

The development of high-resolution shock-capturing schemes has a long history
which we cannot even sketch here (see e.g. the classical references
\cite{Godunov:1959,VanLeer:1979,Harten:1983,Yee:1987}
or the recent textbooks
\cite{LeVeque:1992,GodlewskiRaviart:1996,Kroener:1997,Toro:1999}.
The most common approach is to first develop
a one-dimensional TVD (total-variation-diminuishing) upwind scheme for a scalar
conservation law and then apply it to systems using one-dimensional
characteristic decompositions or approximate Riemann solvers.
Upwind schemes have been used very successfully 
for gas dynamical calculations, where the Riemann problem can be solved
exactly and many approximate Riemann solvers are available. 
For more complicated systems like the granular flow model
considered here characteristic decompositions are often not available,
and the Riemann problem cannot be solved analytically. 
Therefore we have chosen an alternative approach to 
high-resolution shock-capturing, namely the recent non-oscillatory
central (NOC) schemes first introduced by Nessyahu and Tadmor 
\cite{NessyahuTadmor:1990}.
While upwind schemes are higher order extensions of the classical
Godunov scheme, central schemes build upon the (also classical)
Lax-Friedrichs scheme \cite{Lax:1954}. This scheme avoids characteristic decompositions
and Riemann solvers by the use of a staggered grid. When used together
with piecewise constant spatial reconstructions, the Lax-Friedrichs scheme
is more diffusive than Godunov's scheme. However, when one combines the scheme with
TVD-type piecewise linear reconstructions, it becomes competitive with the upwind
schemes. Recently, central schemes have been extended in many directions, see eg.
\cite{ArminjonViallon:1995,ArminjonViallon:1999,JiangTadmor:1998,LieNoelle:2000}
for multidimensional extensions, \cite{RosenbaumRumpfNoelle:2000}
for an adaptive staggered scheme,
\cite{LiuTadmor:1998,LevyPuppoRusso:1998} for third- and higher-order
schemes and \cite{KurganovTadmor:2000,KurganovNoellePetrova:2000}
for central schemes on non-staggered grids, which are precisely at the borderline
of central and upwind schemes.

Here we adapt the second order NOC scheme of Nessyahu and Tadmor to
include an earth pressure coefficient which has a jump discontinuity as the
flow travels from an expanding into a contracting region, and to
treat the source term which is due to the spatially varying
topography and the gravitational force. The resulting scheme 
works well both in smooth regions and at shocks, which are captured 
within two mesh cells and without any oscillations.

Besides the formation of shock fronts in the interior,
avalanches may also have a vacuum front at their margins.
Similarly as for the equations of gas dynamics, the hyperbolic system
degenerates at the vacuum state. Many shock capturing upwind schemes
produce negative heights at these points and subsequently break down
or become completely unstable. While our NOC scheme is remarkably stable
at the margins, it does not capture the vacuum front as well as the Lagrangian
moving mesh scheme.
To overcome this imperfection, we augment the NOC scheme with an
algorithm that tracks the vacuum front. The combined front-tracking
non-oscillatory central scheme is accurate and robust both at
shocks and at the margins of the granular avalanche.

The ensuing analysis commences in \S\,2 with the presentation of the
governing SH-equations in conservative and non-conservative form; then
the jump conditions of mass and momentum at singular surfaces will be
stated and the solution to a single shock wave (a hydraulic
jump) will be presented; \S\,2 closes with the construction of exact
similarity solutions of a parabolic heap moving down a rough incline.
\S\,3 introduces the numerical techniques; at first the Lagrangian
integration technique is described; it is followed by the presentation
of the non-oscillatory central (NOC) scheme.
% the TVD limiters and WENO
%(weighted essentially non-oscillatory) spatial reconstruction.
In \S\,4 we augment the NOC scheme (which uses a fixed Eulerian grid)
with a Lagrangian type front-tracking method in the marginal cells.
\S\,5 elaborates on numerical results. The travelling shock wave cannot be
handled by the Lagrangian method, but the NOC scheme can do so with
very little diffusion accross the shock. On the other hand, the
parabolic similarity solution is well produced by the Lagrangian
integration technique, but much less accurately by the NOC schemes unless
the Lagrangian front-tracking is introduced for the marginal cells. It
is also shown that the NOC scheme with piecewise linear spatial reconstructions
applying standard TVD type slope limiters exhibits some oscillations near smooth
local maxima. We remove these oscillations by incorporating
a piecewise quadratic WENO (weighted essentially non-oscillatory)
reconstruction into our scheme.
Our final numerical experiment combines all the difficulties treated in
the paper: an avalanche with a vacuum front at the margins expands as it
flows downhill and contracts as it hits the flat runout (so the earth
pressure coefficient changes discontinuously at the transition point).
As the avalanche comes to a halt at the bottom, a shock wave develops
and propagates upslope. Our NOC front tracking scheme
handles this challenging flow very satisfactorily.
\S\,6 presents conclusions and gives an outlook to further work.

%
% Section 2
%
\section{Governing Equations}
%\end{center}
%
A detailed derivation of the Savage-Hutter theory has been given in
\cite{SavageHutter:1989,SavageHutter:1991}. Here we confine ourselves to a
brief description.
Although cohesionless granular materials exhibit dilatancy effects
numerous experiments have confirmed that during rapid dense flow
it is reasonable to assume that the avalanche is incompressible
with constant uniform density $\rho_0$."
During flow the body behaves as a
Mohr-Coulomb plastic material at yield.
As the avalanche slides over the rigid basal topography a Coulomb
dry friction force resists the motion. The basal shear stress is
therefore equal to the normal basal pressure multiplied by a coefficient
of friction $\tan\delta$, where $\delta$ is termed the
basal friction angle \cite{KochGreveHutter:1994}.
Scaling analysis isolates the physically significant terms in the
governing equations and identifies those terms that can be neglected.
Plane flow configurations are our focus in this paper, so depth
integration reduces the theory to one spatial dimension.  The leading
order, dimensionless, depth integrated equations for the
local thickness of the avalanche $h$ and the momentum $h u$ 
($u$ is the downslope velocity) reduce to
\begin{align}
\diffx{h}{t}+ \diffx{}{x}(hu)&=0,
\label{2.1}\el
\diffx{(hu)}{t}+\diffx{}{x}(hu^2+\beta_x h^2/2)
&=hs_x
\label{2.2}
\end{align}
with net driving force 
\begin{equation}
\displaystyle
s_x =\sin\zeta -\text{sgn}(u)\tanD
\bigl(\cosZ+\lambda\kappa u^2\bigr)-\varepsilon\cosZ\diffx{z^b}{x},
%\bigl(\cosZ+\lambda\kappa u^2\bigr),
\label{2.3}
\end{equation}
where $x$ is the arc length measured along the avalanche track, $z^b$
denotes the height of the basal topography relative to the track
(usually $z^b=0$ in one spatial dimension) and $\zeta$ and
$\lambda\kappa$ are the local slope inclination angle and
curvature of the track, respectively.
The term sgn$(u)$\ selects the
orientation of the dry Coulomb drag friction, and $\varepsilon
\ll 1$ is the aspect ratio of a typical thickness and length of the
avalanche.  Note that equations (\ref{2.1}) and (\ref{2.2}) are
written in conservative form \cite{GrayWielandHutter:1999}, while in the original
SH theory the smoothness assumption allows the momentum balance
equation to transform to an evolution equation for the velocity, viz.,
\begin{equation}
\Diffx{u}{t}=s_x-\beta_x\diffx{h}{x}-\frac12 h\diffx{\beta_x}{x}.
\label{2.2_1}
\end{equation}
The factor $\beta_x$ is defined as $\beta_x=\varepsilon\cosZ K_x$
 and the earth pressure
coefficient $K_x$ is given by the {\it ad hoc} assumption
\begin{equation}
K_x=
\displaystyle
\left\{
\begin{array}{ll}
\displaystyle
K_{x_{act}}\hskip.2in &\hbox{for}\hskip.2in \partial u/\partial x > 0\,,\el
\displaystyle
K_{x_{pass}}\hskip.2in &\hbox{for}\hskip.2in \partial u/\partial x < 0\,,
\end{array}
\right.
\label{2.4}
\end{equation}
with
\begin{equation}
\displaystyle
K_{x_{act/pass}}
= 2\biggl(1\mp\sqrt{1-\cos^2\phi/\cos^2\delta}\biggr)\sec^2\phi-1\,,
\label{2.5}
\end{equation}
and $\phi$ is the internal friction angle of the granular material.
Note that the values of the earth pressure coefficient $K_x$ are
based on the postulation of a Mohr-Coulomb plastic behaviour for
the cohessionless yield on the basal sliding surface, see Savage
\& Hutter \cite{SavageHutter:1989,SavageHutter:1991}
for details. In this theory the earth pressure
coefficient $K_x$\ is assumed to be function of the velocity gradient,
i.e.  \(K_x=K_x(\partial u/\partial x)\).

The governing equations look like the shallow-water equations, but
because of the jump in the earth pressure coefficients $K_{x_{act/pass}}$,
the source term $s_x$\ and the free boundary at the front and rear
margins, it becomes much more complicated to develop an appropriate
numerical scheme to describe the flow. The original Lagrange
finite-difference scheme \cite{SavageHutter:1989} is implemented 
for the equation system (\ref{2.1}) and (\ref{2.2_1}) in Lagrangian form,
with primitive variables $h$ and $u$.
The shock capturing scheme developed here is applied
to the system in conservative form
(\ref{2.1}) and (\ref{2.2}), where the conserved quantities are
the avalanche thickness $h$\ and the depth integrated momentum $m=hu$.

In vector notation, equations (\ref{2.1}) and (\ref{2.2}) take the form
\begin{equation}
\u w_t+\u f_x=\u s,
\label{2.8}
\end{equation}
where
\begin{equation}
\u w=\left(
\begin{array}{c} 
\displaystyle
h \\
\displaystyle
m 
\end{array}\right),\quad
\u f=\left(
\begin{array}{c} 
\displaystyle
m \\
\displaystyle
m^2/h+\beta_x{h^2}/{2} 
\end{array}\right)\quad\hbox{and}\quad
\u s=\left(
\begin{array}{c} 
\displaystyle
0 \\
\displaystyle
 h\,s_x
\end{array}\right).
\label{2.9}
\end{equation}
This form is more convenient for mathematical analysis than
(\ref{2.1}) and (\ref{2.2}).
%As the shock-capturing scheme is based on a
%stationary uniform grid, under which the velocity $u$\ is undefined
%outside the avalanche body and a discontinuity arise on the front and
%tail margins. However, we cannot avoid the presence of the velocity
%$u$\ in the equations. That's also the main difficulty and problem of
%accuracy for implementation of this theory.  
%
% Subsection 2.1
%
\subsection{Jump Condition and Travelling wave}
\begin{figure}[t!]
  \begin{center}
    \psfrag{h1}[l]{\Large$h^+$}
    \psfrag{h2}[b]{\Large$h^-$}
    \psfrag{u1}[bc]{\Large$u^+$}
    \psfrag{u2}[bc]{\Large$u^-$}
    \psfrag{Vn}[c]{\Large$-V_n$}
    \psfrag{s}[bc]{\Large$\sigma$}
    \psfrag{z}[bc]{\Large$\zeta$}
    \psfrag{ns}[c]{\Large$\V n_\sigma$}
    \includegraphics[angle=0.0,scale=0.45]{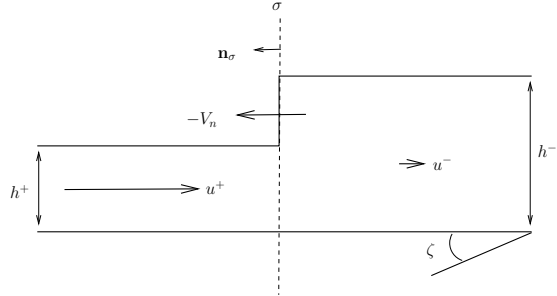}
\vskip12pt
\caption{The plane travelling shock wave can be interpreted as a jump in
  thickness and velocity separating the body of the avalanche into two
  parts on a plane with inclined angle $\zeta$.  ${h^+}$ and ${h^-}$
  are the thicknesses of both sides and $u^+$ and $u^-$ are the
  velocities, respectively, whereas this jump travels with velocity
  $V_n$ up-slopes.
  \label{jumpconditionf}}
  \end{center}
\end{figure}

The Savage-Hutter theory can be used to model the upslope propagating
travelling shock wave observed in experiments
\cite{GrayHutter:1997,GrayTai:1998B} by
introducing the jump conditions (see Fig.\,\ref{jumpconditionf}) of the
balance equations (\ref{2.1}) and (\ref{2.2}) for mass and momentum
\begin{align}
\ljump  h(u-V_n)\rjump &=0\,,\label{jump_1}\\*[8pt]
\TST\ljump hu(u-V_n)+\frac{1}{2}\beta_x h^2\rjump&= 0\,,
\label{jump_2}
\end{align}
where $V_n$ is the normal speed of the singular surface.
Let us suppose that $\ljump\beta_x\rjump=0$ (for example,
this is always satisfied if $\phi=\delta$, i.e.,
$K_{x_{act}}=K_{x_{pass}}$). Substituting (\ref{jump_1}) into
(\ref{jump_2}) (i.e. eliminating $V_n$) yields
the following relation between the depth ratio, $H:=h^-/h^+$, and the
velocity difference
\begin{equation}
(u^+-u^-)^2 = \beta_x h^- \dfrac{H+1}{2} \left(\dfrac{H-1}{H}\right)^2.
\label{VelDiff}
\end{equation}
For an upslope travelling shock wave with travelling wave speed
$V_n$ and corresponding depth ratio $H$, the factor $\beta_x$ is a
function of material and topographic parameters, $\phi$, $\delta$ and
 $\zeta$, which are given by the selected material and topography.
Provided that the depths before and after the shock, $h^+$ and $h^-$,
are known (they can be determined by experiment) and the downslope
velocity is also given (it is normally equal to zero), then the
upslope velocity can be determined by using (\ref{VelDiff})
\begin{equation}
u^+=u^- \; \pm \;  \dfrac{H-1}{H} \left[\beta_x h^- \dfrac{H+1}{2}\right]^{1/2}.
\label{VelDiff2}
\end{equation}
Note that the term under the square root is positive for all positive $H$.
If $H=1$ then $u^+=u^-$, which indicates that no shock wave
(discontinuity) takes place.  Thus, velocity jumps and depth jumps
occur together.

By inspection of the mass balance equation \eqref{jump_1},
the velocity of the shock is given by
\begin{equation}
V_n \; =\; \dfrac{Hu^--u^+}{H-1} \; = \;
u^- \; \mp \;  \left[\beta_x h^- \dfrac{H+1}{2 H^2}\right]^{1/2}.
\label{vel_jump}
\end{equation}
Note that as $h^+$ tends to $h = h^-$, then $u^+$ tends to $u = u^-$ and
\begin{equation}
V_n \; \to \; u  \; \mp \;  \left[\beta_x h\right]^{1/2},
\label{vel_jump2}
\end{equation}
so we have recovered the characteristic speeds of the shallow water equations.
Now we apply Lax' shock inequalities \cite{Lax:1957} to single out the physically
relevant branches of the shock curves: for the first family, with characteristic
speed $u - \sqrt{\beta_x h}$, we require that
$$
u^+ - \left[\beta_x h^+\right]^{1/2} \; > \; V_n \; = \; u^-
- \left[\beta_x h^- \dfrac{H+1}{2 H^2}\right]^{1/2}
\; > \; u^- - \left[\beta_x h^-\right]^{1/2},
$$
which implies $H > 1$ (recall that the upslope state ``+''
lies to the left of the shock). Analogously, for the second family, with
characteristic speed $u + \sqrt{\beta_x h}$, we obtain $H < 1$.
For example, an upward jump ($h^+ < h^-$) can only be carried by a shock
of the first family, and in this case $u^+ > u^- > V_n$, so particles
which cross the shock are condensed and slow down. 

%This travelling shock wave solution serves as a test problem for the
%numerical schemes. It reveals the capability of shock capturing of the
%proposed schemes.
%

% Subsection 2.2
%
\subsection{Similarity Solution}
Consider the motion of a finite mass of granular material along a flat
plane, i.e. $\zeta$ is constant and $\lambda \kappa = 0$
in \eqref{2.3}. In \cite{SavageHutter:1989} one particular similarity solution
to a moving boundary problem of finite mass was derived; this solution
is now generalised (see \cite{Tai:2000}. To this end we introduce a moving coordinate
system with velocity
\begin{equation}
\displaystyle
u_0(t) = u_0(0) + \int^t_0 (\sin\zeta -\tan\delta\,\cos\zeta)\,\d t
\label{2.17}
\end{equation}
on an plane with inclination angle $\zeta$. This velocity is
due to the net driving force $s_x$ in \eqref{2.3}, where we assume that
the velocity is positive for positive times, i.e. sgn$(u) = 1$.
The relative velocity $\breve u$ in the moving
coordinate system is then given by
\begin{equation}
\breve u=u-u_0(t)\point
\label{relVel}
\end{equation}

A symmetric bulk is considered and the origin of the moving coordinate
system is selected to lie at the centre where the surface gradient,
$\dd{h}{x}$, is zero. To keep the symmetric depth profile during the
motion the relative velocity is further assumed to be skew--symmetric,
\(
\breve u(\xi,\,t)=-\breve u(-\xi,\,t),
%\label{symmHU}
\) 
where
\begin{equation}
\xi=x-\int_0^t u_0(t')\d t'
\end{equation}
indicates the distance from the origin in the moving coordinates.
Provided that $g(t)$ is the distance from the coordinate origin to
the margin at time $t$, the physical domain occupied by the granular
mass can be mapped from $[-g(t),\,g(t)]$ to the fixed domain
$[-1,\,1]$ by
\begin{equation}
\eta=\dfrac{1}{g(t)} \left\{ x-\int^t_0 u_0(t')\,\d t'\right\},
\quad\hbox{where}\quad \eta \in[-1,1]. 
\label{7.2}
\end{equation}
With this coordinate mapping, $(x,\,t)~\rightarrow~(\eta,\,\tau)$, the
model equations (\ref{2.1}) and (\ref{2.2}) reduce to
\begin{align}
%\begin{array}{r}
\diffx{h}{t}-\eta\dfrac{g'}{g}\diffx{h}{\eta}+\dfrac{1}{g}\diffx{}{\eta}(h\breve u)&=0\,,
\label{similarsol1}\\*[10pt]
\diffx{\breve u}{t}-\eta\dfrac{g'}{g}\diffx{\breve u}{\eta}
+\dfrac{1}{g}\left(\breve u\diffx{\breve u}{\eta}+\beta_x \diffx{h}{\eta}\right)&=0\,,
%\end{array}
\label{similarsol2}
\end{align}
where the $\tau$ is again replaced by $t$ and we have used
$g'=\d g/\d t = - u_0 / \eta$.  

Now we assume that \( \breve u(\eta, t) \) varies linearly in
$\eta$. Since the margins move with relative speeds $\pm g'(t)$, this
yields \( \breve u(\eta, t)=\eta g'(t) \). Now the evolution
equations (\ref{similarsol1}) and (\ref{similarsol2}) reduce to
\begin{align}
\diffx{h}{t}+\dfrac{g'}{g}h&=0,\label{7.41}\\
\eta g''+\dfrac{\beta_x}{g}\diffx{h}{\eta}&=0,
\label{7.4}
\end{align}
where $g''=\d^2g/\d t^2$. Integrating (\ref{7.4}) subject to the
boundary condition either $h(\eta=1)=0$ or $h(\eta=-1)=0$, it follows
that the thickness is described by
\begin{equation}
h(\eta,\,t)=\dfrac{g(t)g''(t)}{2\beta_x}(1-\eta^2).
\label{parab}
\end{equation}
This implies that the avalanche body keeps a parabolic thickness
distribution during the motion. With the thickness distribution
(\ref{parab}) one can easily obtain the total mass $M$ to be
\begin{equation}
M=\int^{\xi_{Ft}}_{\xi_f}h(\xi,\, t)\,\d\xi
=\int^1_{-1}h(\eta,\, t)g(t)\,\d\eta
=\dfrac{2}{3}\dfrac{g''g^2}{\beta_x}.
\label{totalM}
\end{equation}
Since mass is conserved,
\begin{equation}
0 = \dfrac{d}{dt} M = \dfrac{2 g}{3 \beta_x} (2 g' g '' + g g''').
\end{equation}
This relation can also be derived directly from the mass balance equation
\eqref{7.41}.

Changing the independent variable $t$ to $g(t)$ and letting
$p(t)=g'(t)$, equation (\ref{totalM}) can be written as
\begin{equation}
p\Diffx{p}{g}=\dfrac{K}{g^2},
\label{gfkt}
\end{equation}
where $3\beta_x M=2 K$. The similarity solution is then obtained by
solving (\ref{gfkt}) with initial condition, $g(0)=g_0$ and $p(0)=p_0$
\begin{equation}
p^2(t)
=2K\left(\dfrac{1}{g_0}-\dfrac{1}{g(t)}\right)+p_0^2.
\label{psol}
\end{equation}
With the definition $\alpha_g=\frac{2K}{g_0},~ \beta_g=p_0^2$ and
$G=(\alpha_g + \beta_g) \, g \;$ it follows that
\begin{equation}
\dfrac{\sqrt{G} \, G^{\,\prime}}{\sqrt{G-2K}} = (\alpha_g + \beta_g)^{3/2}.
\label{Gsol}
\end{equation}
We now use the relation
$$
\dfrac{d}{dG} \left[ \sqrt{G}\sqrt{G-2K}+2K \ln (\sqrt{G}+\sqrt{G-2K}) \right]
= \dfrac{\sqrt G}{\sqrt{G-2K}}
$$
and integrate equation (\ref{Gsol}) to yield
\begin{equation}
\begin{array}{l}
\sqrt{G}\sqrt{G-2K}+2K \ln (\sqrt{G}+\sqrt{G-2K})\el
\quad-\left[\sqrt{G}\sqrt{G-2K}+2K \ln (\sqrt{G}+\sqrt{G-2K})\right]_{t=0}
=(\alpha_g+\beta_g)^{3/2}\, t.
\end{array}
\label{gsol}
\end{equation}
With $g_0=1$, $p_0=0$ we obtain the Savage-Hutter solution
\cite{SavageHutter:1989}
\begin{equation}
\sqrt g \sqrt{g-1}+\ln \left(\sqrt g +\sqrt{g-1}\right)
=\sqrt{2K}\, t,
\label{g_24}
\end{equation}
for which $g(t)>1$. Both (\ref{gsol}) and (\ref{g_24}) are
implicit evolution equations for $g(t)$. Once $g(t)$ is deduced, with
the presumption \( \breve u(\eta, t)=\eta g'(t) \), the complete solution
is then given by (\ref{parab}) and (\ref{psol}),
\begin{equation}
\breve u(\eta,\,t)
=\eta\,\left\{2K\left(\dfrac{1}{g_0}-\dfrac{1}{g(t)}\right)+p_0^2\right\}^{1/2},
\qquad h(\eta,\,t)=\dfrac{3M}{4\,g(t)}(1-\eta^2),
\end{equation}
where $\eta$ is defined in \eqref{7.2}.
In the present similarity solution it is presumed that $u/|u|=1$,
which means that $u>0$ for all $t\ge 0$. From (\ref{relVel}) and the
presumption \( \breve u(\eta, t)=\eta g'(t) \) it follows that
\begin{equation}
u(t)=u_0(t)+\breve u(t) >0\quad \Rightarrow\quad g'(t)< u_0(t),
\quad \hbox{for all}\quad t\ge 0.
\label{Vcon}
\end{equation}
It is very important to verify that the velocity is consistent with
condition (\ref{Vcon}) to keep the parabolic similarity solution
valid.  The generalisation (\ref{gsol}) of (\ref{g_24}) was needed to
have exact solutions with non-vanishing initial velocities (for further
details see \cite{Tai:2000}).
%
% Section 3
%
\section{Numerical Scheme}

The numerical schemes employed in this paper are designed to
explicitly solve the system of equations in 1D and we here introduce a
Lagrangian algorithm and an Eulerian shock-capturing NOC
(Non-Oscillatory Central) scheme.

In the Lagrangian technique \cite{SavageHutter:1989,SavageHutter:1991} the
avalanche body is divided into several cells. The purpose is to find
the velocity of the cell boundaries in order to determine the cell
boundary locations for each time step; so it is a moving-grid method,
whereas, the NOC scheme is built on a stationary uniform grid and
gives a high resolution of the shock solutions without any spurious
oscillations near a discontinuity.

In the Lagrangian method the value of the depth $h_j^n$ is defined as
the volume average within the $j^{th}$ cell for time $t^n$, which is
bounded by $b_{j-1}(t)$ and $b_{j}(t)$, and the boundary
$b_j(t)$ moves with the velocity $u_j$. Whilst, in the NOC scheme
the value of the discretised variable $U_j^n,~ U=h,\,m$\ is defined on
the mesh as the volume average within the $j^{th}$ mesh cell centred
at position $x_j$ for time $t^n$, where the $j^{th}$ cell is bounded
by $x_{j+1/2}$\ and $x_{j-1/2}$.
%
% Subsection 3.1
%
\subsection{Lagrangian method}\label{Lscheme}
\begin{figure}[t!]
  \begin{center}
    \includegraphics[angle=0.0,scale=0.53]{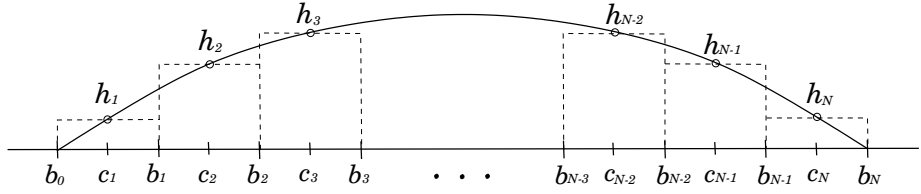}
\vskip12pt
    \caption{The avalanche body is divided into
      $N$ elements with average depth $h_j$, where $c_j$ is the centre
      of the $j^{th}$ element.\label{lagrange}}
\end{center}
\end{figure}

In the Lagrangian method \cite{SavageHutter:1989,SavageHutter:1991} the avalanche
body is divided into $N$ material cells, where $x = b_{j-1}(t)$ and $x = b_j(t)$
denote the boundaries of the cell $j$ at time $t$, see Fig.\,\ref{lagrange}.
These boundaries move with the avalanche velocity, i.e.
$$
\frac{d}{dt} b_j(t) = u_j(t) = u(b_j(t),t).
$$
Integrating the mass balance equation (\ref{2.1}) over the cell yields
\begin{equation}
\begin{array}{l}
\DS
\int\limits_{b_{j-1}}^{b_{j}}\left\{\diffx{h}{t}+\diffx {}x(hu) \right\}\,\d x
= \Diffx{}{t}\int\limits_{b_{j-1}}^{b_{j}}h\, \d x=0\quad
\Rightarrow \quad \Diffx{}{t}V_{cell_j}=0\,,
\end{array}
\label{3.1}
\end{equation}
and implies that the volume (mass) of the cell is conserved during the
motion. Because of this, the mean height of the $j^{th}$ cell can be
determined by
\begin{equation}
\DS h^{n}_j
=\dfrac{V_{cell_j}}{b^{n}_{j}-b^n_{j-1}}.
\label{3.2}
\end{equation}

The computations proceed as follows.  It is assumed that $b^n_j,\ 
h^n_j$ and $u^{n+1/2}_j$ are given as initial values and the new
location of the cell boundary $b^{n+1}_j$ after an elapsed time
 $\triangle t$ is given by
\begin{equation}
b^{n+1}_j=b^{n}_j+\triangle t\; u^{n+1/2}_j.
\label{3.3}
\end{equation}
Note that here the velocity $u_j$ indicates the boundary
velocity of $b_j$.  The momentum balance (\ref{2.2_1}) allows the
velocity of the cell boundary at time $t^{n+1/2}$ to be determined,
\begin{equation}
\hspace{-2mm}
\displaystyle
u^{n+1/2}_j
= u^{n-1/2}_j + \triangle t\left\{s^n_j 
- \varepsilon\cos\!\zeta_j\left(K_x\right)^n_j\!\left(\diffx hx\right)^{\!\!n}_{\!\!j}
-\varepsilon\dfrac{h^n_{j+1/2}}{2}\!
\left(\diffx{(\cos\!\zeta K_x)}{x}\right)^{\!\!n}_{\!\!j}\right\}.
\label{3.4}
\end{equation}
The net driving acceleration $s^n_j$ as given by (\ref{2.3}) is
\begin{equation}
\hspace*{-2mm}
\displaystyle
s^n_j =\sinZ_j -\text{sgn}\!\left(u^{n-1/2}_j\right)\tanD
\left\{\cosZ_j+\lambda\kappa_j {\left(u^{n-1/2}_j\right)}^{\!2}\right\}
-\varepsilon\cos\!\zeta_j\!\left(\diffx{z^b}{x}\right)_{\!j},
\label{3.5}
\end{equation}
where $\zeta_j$ represents the local inclination angle, $\kappa_j$ is
the local curvature, and $z^b$ denotes the local basal topography.
Note that the last term at the right-hand side of (\ref{3.4}) contains the
gradient of the earth pressure coefficient, which is neglected in the
numerical scheme of Savage and Hutter \cite{SavageHutter:1989,SavageHutter:1991}.

The earth pressure coefficient $K_x$ is determined by the {\it ad hoc} definition
\begin{equation}
\DS
\left(K_x\right)^n_j=\left\{\begin{array}{ll}
K_{x_{act}},&\quad\hbox{for}\quad u^{n-1/2}_{j+1}\ \ge\ u^{n-1/2}_j,\\[8pt]
K_{x_{pass}},&\quad\hbox{for}\quad u^{n-1/2}_{j+1}\ <\ u^{n-1/2}_j
\end{array}\right.
\label{LaKwert}
\end{equation}
in \cite{SavageHutter:1989,SavageHutter:1991}. The surface (depth) gradients in
(\ref{3.4}) are determined by the depths of the adjacent elements
\begin{equation}
\DS
\left(\diffx hx\right)^{\!n}_{\!j}
=\dfrac{(h^{n}_{j+1}-h^{n}_{j})}{c^{n}_{j+1}-c^{n}_{j}}
=\dfrac{2(h^{n}_{j+1}-h^{n}_{j})}{b^{n}_{j+1}-b^{n}_{j-1}},
\label{3.6}
\end{equation}
where $c^{n}_{j}$ represents the centre of the $j^{th}$ cell,
$c^{n}_{j}=(b^{n}_{j}+b^{n}_{j-1})/2$, at time $t=t^n$, see
Fig.\,\ref{lagrange}.  The height at the cell boundary, $h_{j+1/2}$,
is given by their mean values in adjacent cells,
%
%\begin{equation}
\( h_{j+1/2}=\frac{1}{2}\left(h_j+h_{j+1}\right),
\)
%\end{equation}
%
and the gradient of the earth pressure coefficient is 
\begin{equation}
\left(\diffx{(\cos\zeta K_x)}{x}\right)^{\!n}_{\!j}
=\dfrac{\cos\zeta_{j+1} \left(K_x\right)^{n}_{j+1}
-\cos\zeta_{j} \left(K_x\right)^{n}_{j}}{c^{n}_{j+1}-c^{n}_{j}}.
\end{equation}

However, while this method is excellent for classical smooth
solutions, it looses numerical stability if shocks develop. Shocks are
initiated when the avalanche velocity is faster than its
characteristic speed and the avalanche front reaches the base of the
slope or a solid wall.  Many detailed investigations about granular
shocks were made by Gray and Hutter \cite{GrayHutter:1997}, in which the shock waves
are concerned to be an important property in the granular flows. To
avoid the numerical instability caused by the shocks, an artificial
viscosity term $\mu\,\partial^2 u/\partial x^2$ is introduced and
added to the right hand side of (\ref{3.5}) for numerical stability,
e.g.  \cite{SavageHutter:1989,SavageHutter:1991} and 
\cite{HutterKochPlussSavage:1995}, where the
artificial viscosity $\mu$ was found to have values between $0.01$ and
$0.03$.

%
%
% Section 3.2
%
\subsection{NOC Scheme}

\begin{figure}[t!]
  \begin{center}
    \includegraphics[angle=0.0,scale=0.48]{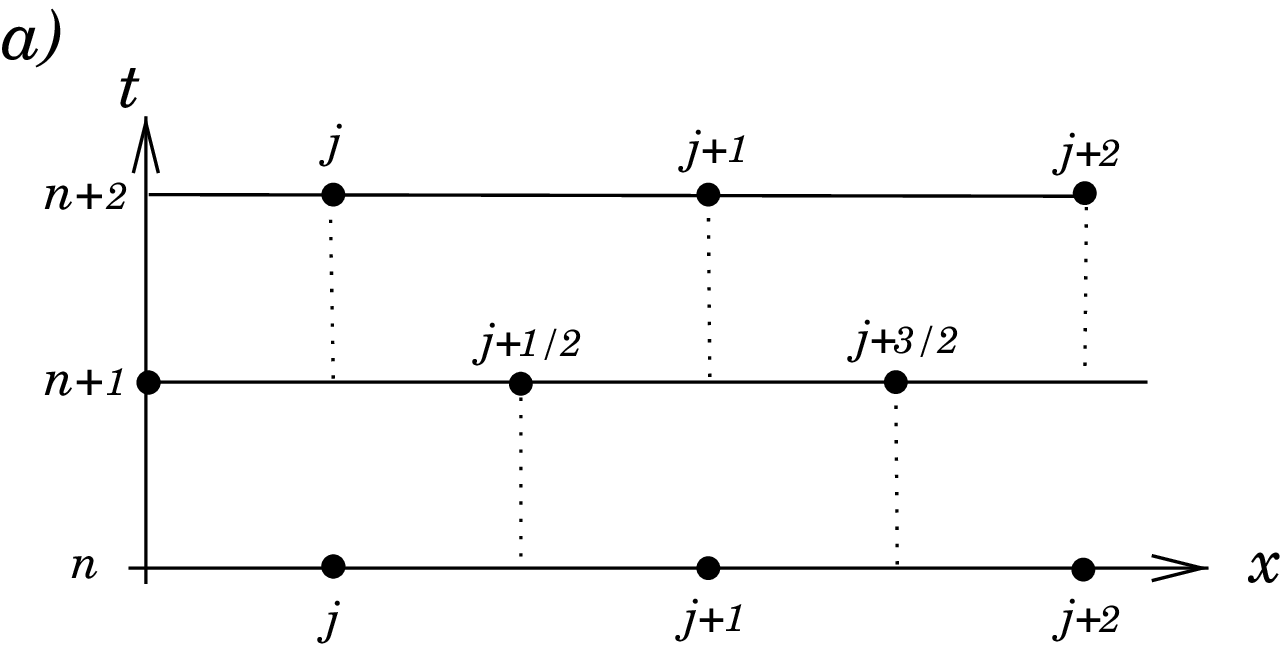}
    \includegraphics[angle=0.0,scale=0.48]{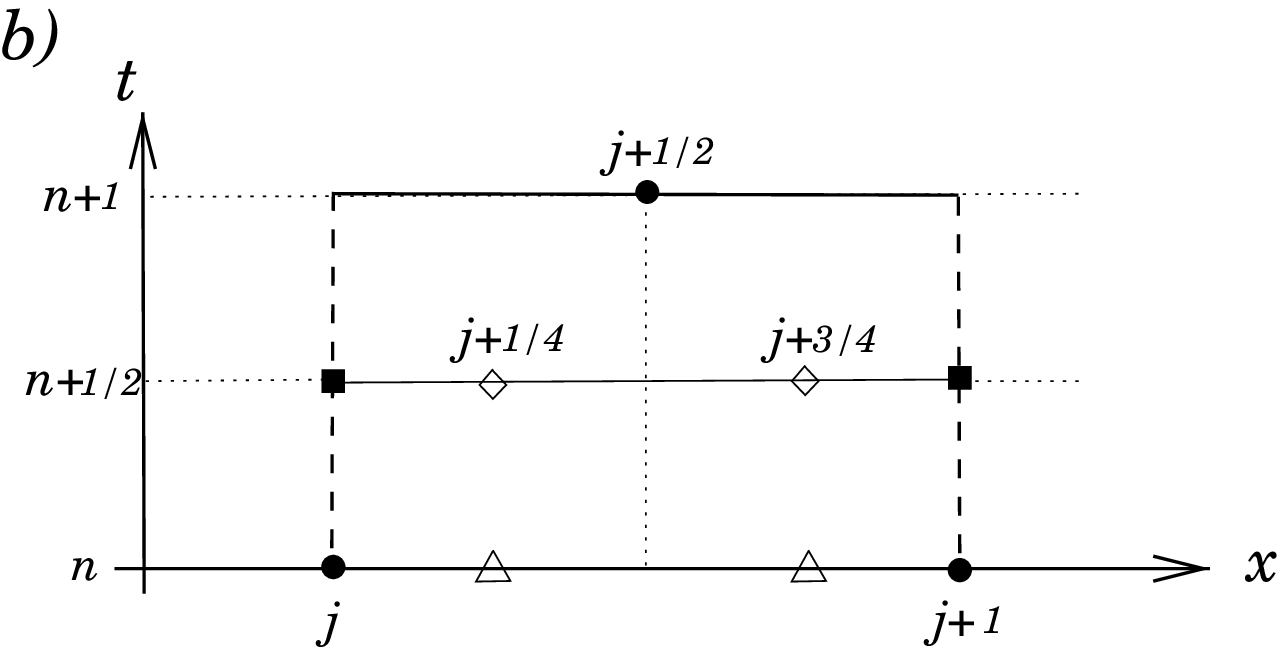}
\vskip12pt
\caption{Diagram of NOCS scheme. 
  a) Grid points computed in the NOCS method.  b) NOCS computational
  diagram, where $\bullet$ indicate the grid points at time level $n$
  and $n+1$, \rule{1.5mm}{1.5mm} represent the quadrature points for the fluxes
  $\V f$ across the cell boundaries, $\diamondsuit$ the quadrature points for
  the source terms $\V s$ and $\triangle$ those for the staggered cell averages
  at the original time $t^n$.
  \label{fig:nocs}}
\end{center}
\end{figure}

The Non-Oscillatory Central Differencing (NOC) scheme  of Nessyahu and
Tadmor \cite{NessyahuTadmor:1990} is a second order accurate extension of the classical
Lax-Friedrichs scheme \cite{Lax:1954}. Let us briefly review the NOC scheme:

We consider the Savage-Hutter equations in the conservative form (\ref{2.8}),
(\ref{2.9}) with $\u w=(h,\,m)^T$ as basic variables.  Let
$\meanW_j^n$ denote the cell average over interval
$[x_{j-\frac12}, x_{j+\frac12}]$ at time $t^n$, and let
\begin{equation}
   \meanW(x, t^n) = \meanW_j^n + \frac{x - x_j}{\Delta x} 
   \meanW'_j
   \label{eq:reconstruction}
\end{equation}
be a piecewise linear reconstruction over the cell,
where $\meanW'_j$ denotes the cell mean derivative determined by a
TVD-limiter \cite{LeVeque:1992} or a central WENO cell reconstruction
\cite{LevyPuppoRusso:1998}.
The main conceptional difference between the NOC schemes and
standard upwind finite difference schemes is the use of a staggered grid.
At time $t^{n+1} = t^n + \Delta t$, the cell averages $\meanW_{j+\frac12}^{n+1}$
are evaluated over the intervals $[x_j, x_{j+1}]$, see Figure \ref{fig:nocs}. 
As a consequence, the boundaries of the cells at the new time level
are the \textit{centers} of the cells at the old time level,
namely the points $x_j$ and $x_{j+1}$.
At these points, the piecewise polynomial reconstruction \eqref{eq:reconstruction}
of the cell averages at the old time level $t^n$ is smooth, and it remains so for
$t < t^{n+1}$ under an appropriate restriction of the timestep (see
\eqref{eq:cfl} below). Therefore, the flux across the boundaries of the cells at the new
time level may be evaluated by Taylor extrapolations using the differential
equation and standard quadrature rules. Here we use the midpoint rule in time
to achieve second order accuracy.  The resulting update takes the form
\begin{equation}
\begin{array}{l}
\meanW_{j+1/2}^{n+1}=\frac{1}{2}\left(\meanW_{j+1/4}^{n}+\meanW_{j+3/4}^{n}\right)
-\frac{\triangle t}{\triangle x}\left(\V f_{j+1}^{n+1/2}-\V f_{j}^{n+1/2}\right)\el
\phantom{\meanW_{j+1/2}^{n+1}=}
+\frac{\triangle t}{2}\left(\V s_{j+1/4}^{n+1/2}+\V s_{j+3/4}^{n+1/2}\right),
\end{array}
\label{3.16}
\end{equation}
as illustrated in Fig.\,\ref{fig:nocs}${}_b$.
The values of $\meanW^n_{j+1/4}$ and $\meanW^n_{j+3/4}$
are determined by the reconstruction \eqref{eq:reconstruction}
over the $j^{th}$ and $(j+1)^{th}$ cell, i.e.
\begin{equation}
%\begin{array}{l}
\TST
\meanW^{n}_{j+ 1/4}=\meanW^n_j     + \frac{1}{4}\meanW'_j,\quad
\meanW^{n}_{j+ 3/4}=\meanW^n_{j+1} - \frac{1}{4}\meanW'_{j+1}.
%\end{array}
\label{3.18}
\end{equation}
The transport flux $\V f$ at the quadrature points $(x_{j},\,t^{n+1/2})$ and
$(x_{j+1},\,t^{n+1/2})$ is approximated by Taylor extrapolation in time,
\begin{equation}
\TST
\V f^{n+1/2}_j=\V f\left(\meanW^{n+1/2}_j\right),
\qquad \meanW^{n+1/2}_{j}=\meanW^n_j 
+ \frac{\triangle t}{2} \left(\dd{\meanW}{t}\right)^{\!n}_{\!j},
\label{fluxNOCS}
\end{equation}
and similarly, the source terms $\V s$ at the quadrature points
$(x_{j+1/4},\,t^{n+1/2})$ and $(x_{j+3/4},\,t^{n+1/2})$
are approximated by space-time Taylor extrapolation
\begin{equation}
\begin{array}{l}
\V s^{n+1/2}_{j+ 1/4}=\V s\left(\meanW^{n+1/2}_{j+ 1/4}\right),\qquad 
\meanW^{n+1/2}_{j+ 1/4}=\meanW^n_{j\phantom{+1}} 
     + \frac{\triangle t}{2} \left(\dd{\meanW}{t}\right)^{\!n}_{\!j\phantom{+1}}
     + \frac{1}{4}\meanW'_j,\\[16pt]
\V s^{n+1/2}_{j+3/4}=s\left(\meanW^{n+1/2}_{j+3/4}\right),\qquad 
\meanW^{n+1/2}_{j+3/4}=\meanW^n_{j+1} 
+ \frac{\triangle t}{2} \left(\dd{\meanW}{t}\right)^{\!n}_{\!j+1}
                        - \frac{1}{4}\meanW'_{j+1}.
\end{array}
\label{sourceNOCS}
\end{equation}
The temporal derivative $(\dd{\meanW}{t})^n_j$ in (\ref{fluxNOCS}) and
(\ref{sourceNOCS}) is determined by using (\ref{2.8}),
\begin{equation}
\left(\dd{\meanW}{t}\right)^{\!n}_{\!j}
=-\left(\dd{\V f}{x}\right)^{n}_{j}+\V s^n_j
=-\V A_j\,\meanW'_j/\dx+\V s^n_j,
\label{3.12}
\end{equation}
where
\begin{equation}
\left(\dd{\V f}{x}\right)^{n}_{j}
=\left(\V A\right)^{n}_{j}\left(\dd{\V w}{x}\right)^{n}_{j}, 
\quad \V A=\dd{\V f}{\V w}
=\left(\begin{array}{cc}
0 & 1 \\
-\frac{m^2}{h^2}+\beta_x h & \frac{2m}{h}
\end{array}\right)
\label{3.13}
\end{equation}
and $\V A$ is the Jacobian of $\V f$. 
Alternatively, one may also use the Jacobian-free approach of
Nessyahu and Tadmor \cite{NessyahuTadmor:1990} and set
$$
\left(\dd{\V f}{x}\right)^{n}_{j} = \V f_{j}' / \Delta x,
$$
where the cell mean derivative $\V f'$ of the flux is again
determined by a TVD-limiter.
Let $a^{max}$ be the maximum wave speed, 
\begin{equation}
a^{max}=\max\limits_{all~j}\left(|u_j|+\sqrt{\beta_{j}h_j}\right),
\quad u_j=m_j/h_j\quad\hbox{for}\quad h_j\neq 0\,.
\end{equation}
The CFL condition
\begin{equation}
\TST
\frac{\dt}{\dx}\left|a^{max}\right| < \frac{1}{2},
\quad\hbox{for all}\quad j
\label{eq:cfl}
\end{equation}
is needed to guarantee that the solution remains smooth at the
space-time quadrature points, so that the Taylor expansions
\eqref{fluxNOCS} and \eqref{sourceNOCS} are justified.

Note that the NOC scheme \eqref{eq:reconstruction} -- \eqref{eq:cfl}
completely avoids the expensive Riemann solvers used in standard upwind schemes
on non-staggered grids. The resulting staggered schemes are easy to code,
computationally efficient and can be applied to general systems of
conservation laws, where the solution of the Riemann problem (i.e. the
initial value problem with piecewise constant data) may complicated or
even impossible.
%
% Section 4
%
\section{Front-tracking method}
\begin{figure}[h!]
\begin{center}
    \includegraphics[angle=0.0,scale=0.8]{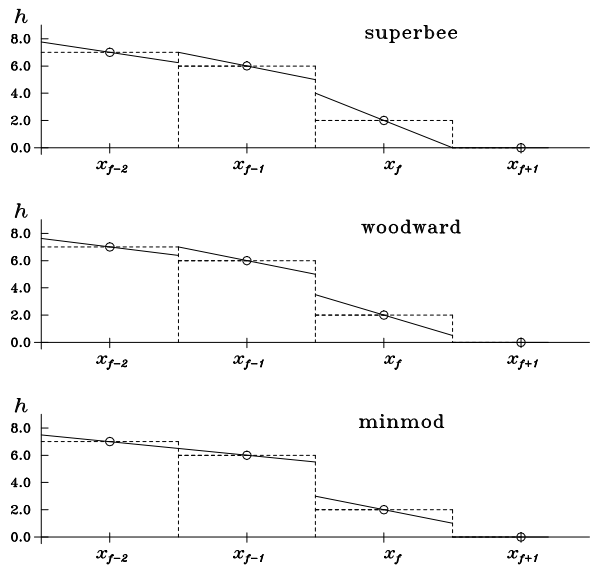}
\vskip12pt
  \caption{
    Example of the depth reconstruction (solid line) determined by
    different TVD limiters, where the circles denote the cell average.
    The front margin lies in the $f^{th}$ cell. In the Eulerian scheme
    one cannot determine where the margin lies. Outside the margin
    there is no material, so that the average depths of the cells
    $f+i,~i\ge1$ are equal to zero.  Different limiters lead to different
    outflows from the avalanche body.
  \label{randcel2}}
\end{center}
%\vskip3mm
\end{figure}

In many applications, the region covered by the granular material has a
finite extension and is limited by a free boundary which moves with the
flow velocity. Outside this region, there is vacuum, so the avalanche
height $h$ and momentum $m$ are zero, and the velocity $u = m / h$ is
not well-defined. The Lagrangian method handles this situation automatically,
since the computational domain moves with the material flow.
The NOC scheme discretizes the differential equations on a stationary
uniform mesh. Note that in general the margin points $x_{Ft}^n$ (the front margin)
and $x_{Tl}^n$ (the tail margin) lie between grid points,
so that it is impossible to point out the margin locations without extra treatment.
Furthermore, it is not straightforward to determine the
proper cell reconstructions over the margin cells.
Fig.\,\ref{randcel2} illustrates an example of depth reconstructions
over the front margin cell determined by various TVD limiters. 
Here and in the following we suppose  that at time $t_n$, the front margin
lies in the $f^{th}$ cell, $x_{f-\frac12} \leq x_{Ft}^n < x_{f+\frac12}$,
and the tail margin in the $t^{th}$ cell,
$x_{t-\frac12} \leq x_{Tl}^n < x_{t+\frac12}$. 

\begin{figure}[b!]
    \psfrag{x}[bl]{\Large $x$}
    \psfrag{h}[bl]{\Large $h$}
    \psfrag{hh}[bl]{\large $\tilde h_f(x)$}
    \psfrag{h1}[bl]{\large $h_2(x)$}
    \psfrag{h2}[bl]{\large $h_1(x)$}
    \psfrag{xf}[bl]{\large $x_{Ft}$}
    \psfrag{r}[bl]{\large $x_f$}
    \psfrag{text}[bc]{\large $\tilde h_f(x)=h_1(x)$ or $h_2(x)$ ?}
    \psfrag{r-1}[bl]{\large $x_{f-1}$}
    \begin{center}
    \includegraphics[angle=0.0,scale=0.6]{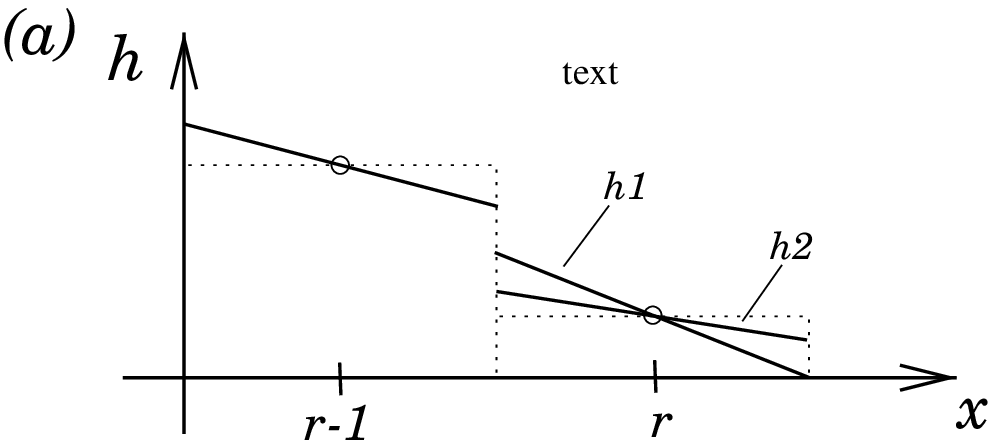}
    \includegraphics[angle=0.0,scale=0.6]{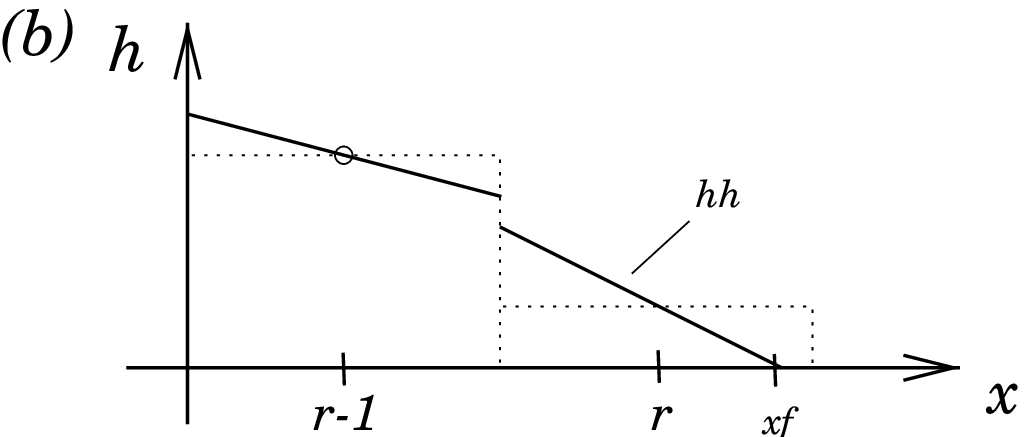}
\vskip12pt
  \caption{
    The reconstruction of the depth $\tilde h_f(x)$ within the margin
    $f^{th}$ cell.  (a) Cell reconstructions based on TVD-limiters
    cannot determine the location of the margin point. Non-appropriate
    reconstructions over the margin cell may result in wrong
    values of the flux at the gridpoint $x_f$, which may
    cause too much outflow from the avalanche body.
    (b) Our front tracking method uses the unique piecewise linear
    reconstruction $\tilde h_f(x)$ over the margin cell which vanishes at the 
    margin point $x_{Ft}$ and preserves the cell average.
    Thus, a reasonable flux at $x_f$ is expected.
  \label{randcell}}\end{center}
\end{figure}

Since our quadrature rule for the fluxes \eqref{fluxNOCS},
\eqref{3.12}, \eqref{3.13} uses a Taylor expansion of the solution,
different limiters will lead to different values of the integrals of
the fluxes across $x_f \times [t_n, t_{n+1}]$ and $x_t \times [t_n, t_{n+1}]$.
To complicate the situation even further, part of these boundaries
may lie in the vacuum region. Note that the fluxes across these
boundaries determine the outflow from the avalanche body,
so non-appropriate cell reconstructions over the margin cell may cause too much
outflow from the avalanche body or even result in a negative depth
around the margin, see Fig.\,\ref{randcell}${}_a$. Thus, the
difficulty is not only to determine the correct numerical flux at the
grid point $x_f$, the wrong numerical flux may also cause
vast stability problems. Adding a thin layer over the whole
computational domain can circumvent the numerical stability problem,
but it is then difficult to determine the locations of the margins,
and the numerical flux out of the avalanche body would even become
unexpectedly large, which results in large numerical diffusion, while
there will be permanent outflow from the avalanche body.  Therefore, a
more refined treatment is needed for the evolution of the 
avalanche margins.

In \cite{Munz:1994}, Munz developed a method to track vacuum fronts in
gas dynamics. His approach is based on appropriate reconstructions
of cell averages behind the front, and the solution of a vacuum
Riemann problem which is used to track the margin locations
at every time step. Here we develop an alternative front tracking
method, which is based on a piecewise linear spatial
reconstruction of the conservative variables up to the front
and Taylor extrapolations in time. Contrary to \cite{Munz:1994}
our approach is Riemann-solver free and therefore fits perfectly
into the framework of central schemes.

The structure of our front-tracking algorithm is as follows:
At the beginning of each timestep (at time $t_n$), the cell averages
$\meanW_j^n$ of the conservative variables and the position of the
margin points $\xft^n$ (front) and $\xtl^n$ (tail) are given.
In the first step, a piecewise linear reconstruction of the data is defined,
the front (tail) velocity is determined and the front (tail)
is propagated from time $t_n$ to $t_{n+1}$.
In the second step, the conservative variables are updated via
\begin{equation}
\hspace*{-5mm}
\begin{array}{l}
\DS \meanW^{n+1}_{j-1/2}=
\dfrac{1}{\triangle x}\int^{x_{j}}_{x_{j-1}}\hspace{-2mm}\V w(x,\,t^n)\,\d x
%+\dfrac{1}{\triangle x}\int^{x^n_{j}}_{x_{j-1/2}}\hspace{-2mm}\V w(x,\,t^n)\,\d x
-\dfrac{1}{\triangle x}
\int^{t^{n+1}}_{t^n}\hspace{-6mm}\left\{\V f(x_{j},t)-\V f(x_{j-1},t)\right\}\,\d t
\\*[16pt]\DS\phantom{h^{n+1}_{j-1/2}=}
+\dfrac{1}{\triangle x}\int^{t^{n+1}}_{t^n}\hspace{-3mm}
\int^{x_{j}}_{x_{j-1}}\hspace{-3mm}\V s(x,t)\,\d x\,\d t,
\end{array}
\label{wf-0}
\end{equation}
Away from the front, the integrals are evaluated by the midpoint rule as in
(\ref{3.16}). Special care has to be taken in the two margin cells (the cells
containing the front and the tail). Each of the integrals on the RHS
of (\ref{wf-0}) may contain parts of the vacuum region. Therefore, we need to replace the
midpoint rule by more delicate quadrature rules over the region covered by the
granular material.    

In order to guide the reader through the details of the algorithm, we would like to
give an outline of the rest of this section. In Section \ref{sect:4.1},
a particular piecewise linear reconstruction of the conservative variables
near the front is derived. In Section \ref{sect:4.2}, the front velocities are computed,
and the fronts are propagated to the new time level. In Section \ref{sect:4.3},
four cases are distinguished for the location of the front relative to
the fixed underlying grid, and their geometry is discussed.
In Sections \ref{sect:4.4}, \ref{sect:4.5} and \ref{sect:4.6}, the three
integrals on the RHS of (\ref{wf-0}) are treated: the data, the fluxes and the source terms.
In Section \ref{sect:4.7}, a special space-time Taylor extrapolation of the conservative
variables near the front is derived, which is needed to compute the solution at the
space-time quadrature points of the three integrals.
Section \ref{sect:4.8} summarizes the algorithm. 
%
%A first reading might focus on Sections \ref{sect:4.1}, \ref{sect:4.2}, \ref{sect:4.3}
%and \ref{sect:4.7} and just browse through the details of Sections \ref{sect:4.4} --
%\ref{sect:4.6} and \ref{sect:4.8}.

%
%

%
% Section 4.1
%
\subsection{Reconstructing the conservative variables}
\label{sect:4.1}
Suppose as before that the front margin is contained in the $f^{th}$ cell,
and the rear margin in the $t^{th}$ cell,
$$
x_{Ft}^n \in [x_{f-\frac12},x_{f+\frac12}],
$$
$$
x_{Tl}^n \in [x_{t-\frac12},x_{t+\frac12}].
$$
We require that the piecewise linear reconstruction $\meanW(x,t_n)$
satisfies the following two criteria:
\begin{itemize}
\item first, it should vanish at the margin points, and
\item second, it should preserve the cell averages.
\end{itemize}
These criteria uniquely determine the reconstruction in the margin cells.
If we denote the cell averaged depths of the front margin and
rear margin cells by $h_f$ and $h_t$,
the depth reconstruction is defined by
\begin{equation}
\tilde h_f(x) = \sigma^h_{f} (x-x_{Ft});
\qquad \sigma^h_{f} 
= \dfrac{-2\,h_f\,\dx}{(x_{Ft}-x_{f-1/2})^2},\quad \hbox{for}~x\in [x_{f-1/2},\,x_{Ft}],
\label{6.11}
\end{equation}
for the front margin cell and by
\begin{equation}
\tilde h_t(x) = \sigma^h_{t} (x-x_{Tl});
\qquad \sigma^h_{t} 
= \dfrac{2\,h_t\,\dx}{(x_{t+1/2}-x_{Tl})^2},\quad \hbox{for}~x\in [x_{Tl},\,x_{t+1/2}],
\label{6.12}
\end{equation}
for the rear margin cell. Outside the margin the depth is equal to
zero. The $x_{Ft}$ and $x_{f-1/2}$ represent the locations of the front
and the internal boundary of the front margin cell, respectively (see
Fig.\,\ref{randcell}).  The $x_{Tl}$ and $x_{t+1/2}$ denote the
locations of the rear and the internal boundary of the rear margin
cell, respectively. The reconstruction of $m=hu$, $\widetilde m(x)$,
is defined analogously.
%
% Section 4.2
%
\subsection{Propagating the front}
\label{sect:4.2}
Besides being very natural and simple, our definition of the reconstructions
over the margin cells has the advantage of leading to uniquely defined, constant values
of the velocities over the whole margin cells,
$$
u(x,t_n) = \dfrac{\widetilde m_{f/t}(x)}{\widetilde h_{f/t}(x)}
= \dfrac{\sigma^m_{f/t}(x-x_{Ft/Tl})}{\sigma^h_{f/t}(x-x_{Ft/Tl})}
= \dfrac{m_{f/t}}{h_{f/t}}.
$$
It is therefore natural to
define the margin velocities by these constant values,
\begin{equation}
\DS
u_{Ft/Tl}^n = \dfrac{m_{f/t}}{h_{f/t}}.
\label{eq:u_FT_n}
\end{equation}
Alternatively, we may also approximate the margin velocity at time
$t^{n+1/2}$ in order to evolve the margin with second order accuracy.
Using the evolution equation \eqref{2.2_1} for the velocity we define
\begin{equation}
u_{Ft/Tl}^{n+1/2} = u_{Ft/Tl}^n + \frac{\triangle t} 2
\left(s_{f/t}^n-\beta_x \sigma^h_{f/t} \right).
\label{eq:u_FT_n+half}
\end{equation}
Here we have used the fact that $h$ vanishes at the front.
The location of the margin at the new time level is then given by
\begin{equation}
x^{n+1}_{Ft/Tl}=x^{n}_{Ft/Tl}+\triangle t\;u_{Ft/Tl},
\label{6.3}
\end{equation}
where either $u_{Ft/Tl} = u_{Ft/Tl}^n$ or
$u_{Ft/Tl} = u_{Ft/Tl}^{n+1/2}$.
%
% Section 4.3
%
\subsection{Intersecting the front and the grid}
\label{sect:4.3}
\begin{figure}[h]
    \begin{center}
    %\begin{minipage}{3cm} 
    \psfrag{(b)}[bl]{\large }
    \psfrag{F}[bl]{\Large $\V s^{n+1/2}_{f-1}$}
    \psfrag{F2}[bl]{\Large $\V s^{n+1/2}_{f}$}
    \psfrag{fn}[bl]{\Large $x^{n+1}_{Ft}$}
    \psfrag{fx}[bl]{\Large $x^{n}_{Ft}$}
    \psfrag{x}[bl]{\Large $x$}
    \psfrag{t}[bl]{\Large $t$}
    \psfrag{n+1}[bl]{\Large $t^{n+1}$}
    \psfrag{n}[bl]{\Large $t^n$}
    \psfrag{r12}[bc]{\Large $x_{f-1/2}$}
    \psfrag{r1}[bl]{\Large $x_{f-1}$}
    \psfrag{r}[bl]{\Large $x_{f}$}
    \psfrag{rp}[bl]{\Large $x_{f+1}$}
    \psfrag{rp2}[bc]{\Large $x_{f+1/2}$}
    \psfrag{sa}[bl]{\Large $\V s_{Ft}(t)=0$}
     \includegraphics[angle=0.0,height=4cm,width=12cm]{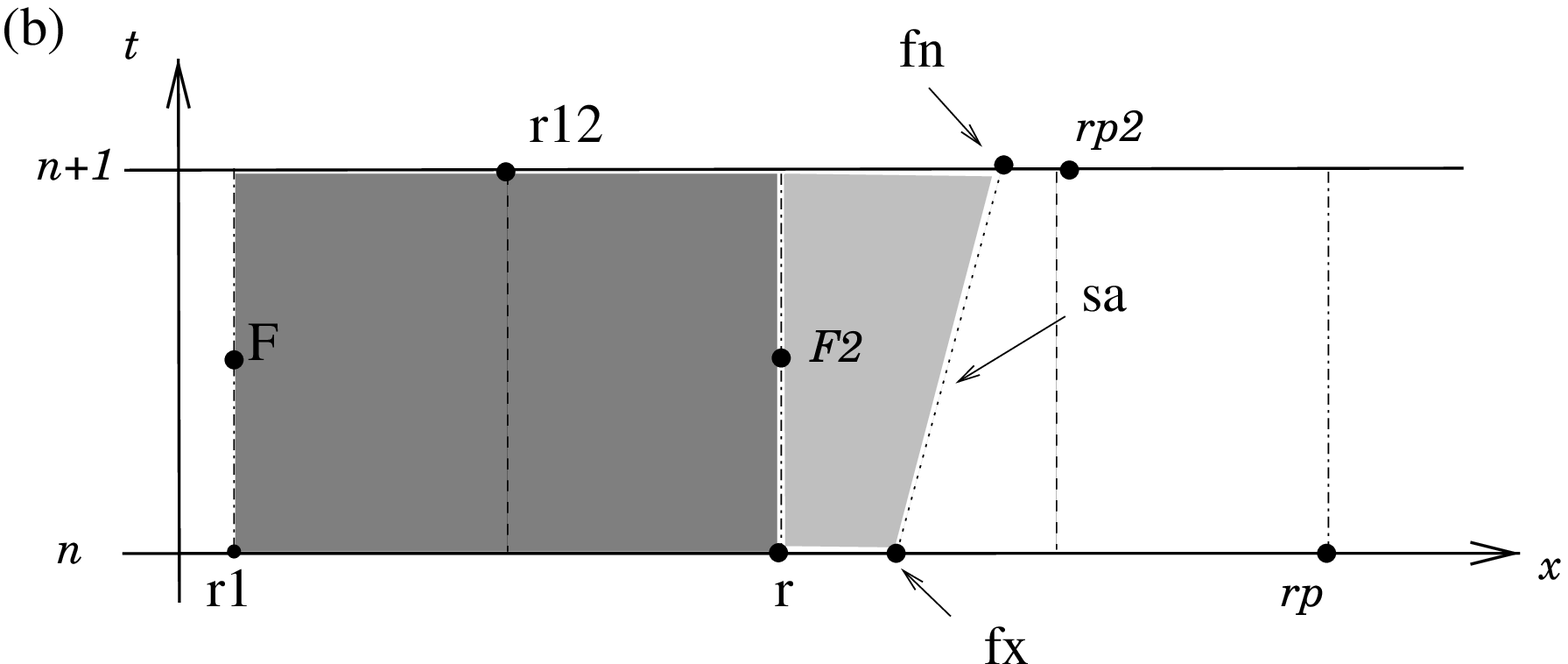}

     case I

     %\end{center}
     %\end{minipage}
     %\hfill
     %\begin{minipage}{3cm}

    \vspace*{1ex}

    \psfrag{(b)}[bl]{\large }
    \psfrag{F}[bl]{\Large $\V s^{n+1/2}_{f-1}$}
    \psfrag{F2}[bl]{\Large $\V s^{n+1/2}_{f}$}
    \psfrag{fn}[bl]{\Large $x^{n+1}_{Ft}$}
    \psfrag{fx}[bl]{\Large $x^{n}_{Ft}$}
    \psfrag{x}[bl]{\Large $x$}
    \psfrag{t}[bl]{\Large $t$}
    \psfrag{n+1}[bl]{\Large $t^{n+1}$}
    \psfrag{n}[bl]{\Large $t^n$}
    \psfrag{r12}[bc]{\Large $x_{f-1/2}$}
    \psfrag{r1}[bl]{\Large $x_{f-1}$}
    \psfrag{r}[bl]{\Large $x_{f}$}
    \psfrag{rp}[bl]{\Large $x_{f+1}$}
    \psfrag{rp2}[bc]{\Large $x_{f+1/2}$}
    \psfrag{sa}[bl]{\Large $\V s_{Ft}(t)=0$}
    %\begin{center}
    \includegraphics[angle=0.0,height=4cm,width=12cm]{figures/work/gibbs/tai/Diss/image/FBsourceb.eps}

     case II
     
     %\end{center}
     %\end{minipage}

    \vspace*{1ex}

    \psfrag{(b)}[bl]{\large }
    \psfrag{F}[bl]{\Large $\V s^{n+1/2}_{f-1}$}
    \psfrag{F2}[bl]{\Large $\V s^{n+1/2}_{f}$}
    \psfrag{fn}[bl]{\Large $x^{n+1}_{Ft}$}
    \psfrag{fx}[bl]{\Large $x^{n}_{Ft}$}
    \psfrag{x}[bl]{\Large $x$}
    \psfrag{t}[bl]{\Large $t$}
    \psfrag{n+1}[bl]{\Large $t^{n+1}$}
    \psfrag{n}[bl]{\Large $t^n$}
    \psfrag{r12}[bc]{\Large $x_{f-1/2}$}
    \psfrag{r1}[bl]{\Large $x_{f-1}$}
    \psfrag{r}[bl]{\Large $x_{f}$}
    \psfrag{rp}[bl]{\Large $x_{f+1}$}
    \psfrag{rp2}[bc]{\Large $x_{f+1/2}$}
    \psfrag{sa}[bl]{\Large $\V s_{Ft}(t)=0$}
    %\begin{center}  
    \includegraphics[angle=0.0,height=4cm,width=12cm]{figures/work/gibbs/tai/Diss/image/FBsourceb.eps}

     case III

    \vspace*{1ex}

    \psfrag{(b)}[bl]{\large }
    \psfrag{F}[bl]{\Large $\V s^{n+1/2}_{f-1}$}
    \psfrag{F2}[bl]{\Large $\V s^{n+1/2}_{f}$}
    \psfrag{fn}[bl]{\Large $x^{n+1}_{Ft}$}
    \psfrag{fx}[bl]{\Large $x^{n}_{Ft}$}
    \psfrag{x}[bl]{\Large $x$}
    \psfrag{t}[bl]{\Large $t$}
    \psfrag{n+1}[bl]{\Large $t^{n+1}$}
    \psfrag{n}[bl]{\Large $t^n$}
    \psfrag{r12}[bc]{\Large $x_{f-1/2}$}
    \psfrag{r1}[bl]{\Large $x_{f-1}$}
    \psfrag{r}[bl]{\Large $x_{f}$}
    \psfrag{rp}[bl]{\Large $x_{f+1}$}
    \psfrag{rp2}[bc]{\Large $x_{f+1/2}$}
    \psfrag{sa}[bl]{\Large $\V s_{Ft}(t)=0$}
  
    \includegraphics[angle=0.0,height=3cm,height=4cm,width=12cm]{figures/work/gibbs/tai/Diss/image/FBsourceb.eps}

     case IV

   \vskip12pt

  \caption{
    The four cases for the propagation of the front margin. \\[2ex]
    \textbf{\large In cases I, III and IV replace by correct figures!! }                                       
    \label{fig:regions}}
  \end{center}
\end{figure}

Once the new location of the margin is given, the new margin cell at
the next time step is then determined. The CFL condition
\eqref{eq:cfl} guarantees that
\( \left| u^{n}_{Ft/Tl}\,\triangle t \right| < \triangle x /2 \),
so the margin point $x_{Ft/Tl}$ can at most pass through gridpoint $x_{f/t}$
during one time step.  For example, with this
condition the front can only lie in either one of the two adjacent
cells of the margin $f^{th}$ cell, which are the $(f-\halb)^{th}$ and
y$(f+\halb)^{th}$ cells, see Fig.\,\ref{fig:regions}.  There are four
possible cases for the motion of the front margin point:
\begin{itemize}
\item \textbf{case I : } $x^{n}_{Ft} \leq x_f$, and $x_{f-1} <  x^{n+1}_{Ft}\leq x_{f}$,
\item \textbf{case II :} $x^{n}_{Ft} > x_f$, and $x_f<x^{n+1}_{Ft}\leq x_{f+1}$,
\item \textbf{case III : } $x^{n}_{Ft} \leq x_f$,  and $x_f<x^{n+1}_{Ft}\leq x_{f+1/2}$,
\item \textbf{case IV : } $x^{n}_{Ft}> x_f$, and $x_{f-1/2}< x^{n+1}_{Ft}\leq x_{f}$,
\end{itemize}
where $x^{n}_{Ft}$ and $x^{n+1}_{Ft}$ are the front locations at
$t^{n}$ and $t^{n+1}$, respectively. In cases I and II, the front does not pass gridpoint $x_f$,
while in the cases III and IV it does, see Figure \ref{fig:regions}. 
%
% Fig.
%
In each case we have to determine the cell averages of the relevant cells
$\meanW^{n+1}_{f-1/2}$ and $\meanW^{n+1}_{f+1/2}$ by
integrating the governing equations over
$[x_{f-1},x_{f}]\times[t^n,t^{n+1}]$ and
$[x_{f},x_{f+1}]\times[t^n,t^{n+1}]$, respectively,
i.e. we have to evaluate the three integrals on the RHS of (\ref{wf-0}).
These integrals involve the data $\V w$, the fluxes $\V f$ and the
source term $\V s$. In the following, we derive quadrature rules
which are exact for linear functions. The tail margin can be treated
completely analogously.
%
% Section 4.4
%
\subsection{The integral of the data}
\label{sect:4.4}
First we integrate the linear reconstruction $\V w(x,t_n)$ of the data at time
$t_n$ over the interval $[x_{f-1}, x_{f}]$. In cases I and III, this interval contains
the front, while it does not in cases II and IV. We obtain
\begin{equation}
\dfrac{1}{\triangle x}\int^{x_{f-1}}_{x_{f}}\hspace{-2mm}\V w(x,\,t^n)\,\d x =
\left\{
\begin{array}{cl}
\half \meanW_{f-3/4} + \meanW_f^n & \textrm{in cases I and III} \\[1ex]
\half (\meanW_{f-3/4} + \meanW_{f-1/4}^n)  & \textrm{in cases II and IV}
\end{array}
\right.
\label{eq:4.4.1}
\end{equation}
Here $\meanW_{f-3/4}$ is given by (\ref{3.18}), and
$\V w^n_{f-1/4}$ by (\ref{6.11}). With a given
front location $x^n_{Ft}$ and $\meanW^n_{f}$ it is
\begin{equation}
\V w^n_{f-1/4}
=2\,\meanW^n_f\left\{
1-\left(\frac{x^n_{Ft}-x_f}{x^n_{Ft}-x_{f-1/2}}\right)^2\right\}.
\label{eq:4.4.2}
\end{equation}
Next we consider the integral over the interval $[x_{f}, x_{f+1}]$.
Using (\ref{6.11}) once more we obtain
\begin{equation}
\dfrac{1}{\triangle x}\int^{x_{f}}_{x_{f+1}}\hspace{-2mm}\V w(x,\,t^n)\,\d x =
\left\{
\begin{array}{cl}
0 & \textrm{in cases I and III} \\[1ex]
\meanW^n_f \left(\frac{x^n_{Ft}-x_f}{x^n_{Ft}-x_{f-1/2}}\right)^2  & \textrm{in cases II and IV}
\end{array}
\right.
\label{eq:4.4.3}
\end{equation}
%
%
% Section 4.5
%
\subsection{The integral of the fluxes}
\label{sect:4.5}
Due to the restriction of the timestep, the only grid-position which is
possibly intersected by the front during the time-interval $[t_n,t_{n+1}]$  
is $x = x_f$. Therefore, the flux at $x_{f-1}$ can be evaluated exactly as 
in the interior of the domain,
\begin{equation}
\dfrac{1}{\triangle t}\int^{t_{n}}_{t_{n+1}}\hspace{-2mm}\V f (\V w(x_{f-1},\,t))\,\d t =
\V f_{f-1}^{n+1/2},
\label{eq:4.5.1}
\end{equation}
where $f_{f-1}^{n+1/2}$ is given by (\ref{fluxNOCS}). The flux at $x_{f+1}$ vanishes,
since this point lies in the vacuum region during the whole time interval.
It remains to compute the flux at $x_f$. In cases III and IV, where the front 
crosses $x_f$, we use the midpoint-rule in time over
that part of the  interface which lies within the
region covered by granular material. Let $\bar t$ and $\overline{\triangle t}$ be the midpoint
and the length of this time interval.
If $t^*$ is the time at which the front intersects $x_f$,
defined by
\begin{equation}
x_{Ft}^n + (t^* - t_n) u_{Ft}^n = x_f, 
\label{eq:4.5.2}
\end{equation}
then
\begin{equation}
\bar t =
\left\{
\begin{array}{cl}
(t_{n+1} + t^*)/2 & \textrm{in case III} \\[1ex]
(t_{n} + t^*)/2  & \textrm{in case IV}
\end{array}
\right.
\label{eq:4.5.3}
\end{equation}
and 
\begin{equation}
\overline{\triangle t} =
\left\{
\begin{array}{cl}
t_{n+1} - t^* & \textrm{in case III} \\[1ex]
t^* - t_n  & \textrm{in case IV}
\end{array}
\right.
\label{eq:4.5.4}
\end{equation}
The midpoint rule for the flux now gives
\begin{equation}
\dfrac{1}{\triangle t}\int^{t_{n}}_{t_{n+1}}\hspace{-2mm}\V f (\V w(x_{f},\,t)\,\d t =
\left\{
\begin{array}{cl}
0 & \textrm{in case I} \\[1ex]
\V f_f^{n+1/2} & \textrm{in case II} \\[1ex]
\frac{\overline{\triangle t}}{\triangle t} \V f_f^{\bar t} & \textrm{in cases III and IV}
\end{array}
\right.
\label{eq:4.5.5}
\end{equation}
Here $\V f_f^{\bar t} = \V f(\V w(x_f, \bar t))$.
In Section \ref{sect:4.7} below we will extrapolate the solution $\V w$ to the quadrature
point $(x_f, \bar t)$.
%
%
% Section 4.6
%
\subsection{The integral of the source term}
\label{sect:4.6}
The source term $\V s$ has to be integrated over the quadrilateral regions
shown in Figure \ref{fig:regions}. Let us call these areas of integration 
$\Omega$. In the following lemma, we give a quadrature rule 
which is exact for linear functions vanishing at the front.
\begin{lemma}
\label{lemma:4.1}
Let $a,b,\tau \geq 0$ and 
$$
\Omega := \{(x,t): \hat t \leq t \leq \hat t + \tau,
\hat x \leq x \leq \hat x + a + (b-a)(t-\hat t)/\tau \}.
$$
Let $\V s$ be a linear function over $\Omega$ which vanishes at the boundary
$x = \hat x + a + (b-a)(t-\hat t)/\tau$. Then
\begin{equation}
\int\int\limits_\Omega \V s(x,t) dx dt =
\frac13 \tau \frac{a^2 + a b + b^2}{a+b} \V s(\hat x, \hat t + \tau/2)
=: \omega \V s(\hat x, \hat t + \tau/2).
\label{eq:4.6.1}
\end{equation}
\end{lemma}
\textbf{Proof:}
W.l.o.g. let $\hat x = \hat t = 0$. The general form of $\V s$ is given by
$$
\V s(x,t) = (x - a - (b-a) t/\tau) \sigma
$$
where $\sigma$ is a real constant. W.l.o.g. let $\sigma = 1$.
Then a direct computation gives that
$$
\int\int\limits_\Omega \V s(x,t) dx dt = - \frac\tau6(a^2+ab+b^2) = \omega \V s(0,\frac\tau2)
$$
$\square$

Equation (\ref{eq:4.6.1}) may be
interpreted as a special quadrature rule with node $(\hat x, \hat t + \tau/2)$.
We have choosen this node because it appears also in the quadrature rule for the fluxes
treated in Section \ref{sect:4.5} above, so we can minimize the evaluations of the
solution $\V w$.
 
In the following we apply the lemma to the four cases.
Let $\overline\Omega$ be the region covered by the granular material.
First, we compute the integral over the intersection of $\overline\Omega$
with the union of the $(f-1/2)^{th}$ and the $(f+1/2)^{th}$ cell,
$\Omega = \overline\Omega \cap ([x_{f-1},x_{f+1}]\times[t^n,t^{n+1}])$.
Using $\hat x = x_{f-1}$, $\hat t = t_n$,
$a = x_{Ft}^n - x_{f-1}$, $b = x_{Ft}^{n+1} - x_{f-1}$ and $\tau = \triangle t$
in Lemma \ref{lemma:4.1} gives
\begin{equation}
\hspace*{-5mm}
\int^{t^{n+1}}_{t^n}\hspace{-3mm}
\int^{x_{f+1}}_{x_{f-1}}\hspace{-3mm}\V s(x,t)\,\d x\,\d t
= \omega_{f-1} \V s_{f-1}^{n+1/2} 
\label{eq:4.6.2}
\end{equation}
with
\begin{equation}
\hspace*{-5mm}
\omega_{f-1} = \frac{\triangle t}{3} \frac{(x_{Ft}^{n}-x_{f-1})^2 +
(x_{Ft}^{n}-x_{f-1})(x_{Ft}^{n+1}-x_{f-1}) + (x_{Ft}^{n+1}-x_{f-1})^2}
{x_{Ft}^{n} + x_{Ft}^{n+1}- 2 x_{f-1}}.
\label{eq:4.6.3}
\end{equation}
Similarly, for the integral over $\overline\Omega \cap 
([x_{f},x_{f+1}]\times[t^n,t^{n+1}])$ we obtain
\begin{equation}
\hspace*{-5mm}
\int^{t^{n+1}}_{t^n}\hspace{-3mm}
\int^{x_{f+1}}_{x_{f}}\hspace{-3mm}\V s(x,t)\,\d x\,\d t =
\omega_f \V s_f^{\bar t},
\label{eq:4.6.4}
\end{equation}
where $\bar t = t^{n+1/2}$ in cases I and II and $\bar t$ is given by (\ref{eq:4.5.3})
in cases III and IV, and the wheight is given by 
\begin{equation}
\hspace*{-5mm}
\omega_f =
\left\{
\begin{array}{cl}
0 & \textrm{in case I}\\[2ex]
\frac{\triangle t}{3} \frac{(x_{Ft}^{n}-x_{f})^2 + (x_{Ft}^{n}-x_{f})(x_{Ft}^{n+1}-x_{f})
+ (x_{Ft}^{n+1}-x_{f})^2} {x_{Ft}^{n} + x_{Ft}^{n+1}- 2 x_{f}} 
  & \textrm{in case II} \\[2ex]
\frac{\overline{\triangle t}}{3}(x_{Ft}^{n+1}-x_{f}) & \textrm{in case III} \\[1ex]
\frac{\overline{\triangle t}}{3}(x_{Ft}^{n}-x_{f})   & \textrm{in case IV}
\end{array}
\right.
\label{eq:4.6.5}
\end{equation}
Here  $\overline{\triangle t}$ is given by (\ref{eq:4.5.4}).
The integral over $[x_{f-1},x_{f}]\times[t^n,t^{n+1}]$ is then computed
by substracting the integral over $[x_{f},x_{f+1}]\times[t^n,t^{n+1}]$
from that over $[x_{f-1},x_{f+1}]\times[t^n,t^{n+1}]$,
\begin{equation}
\hspace*{-5mm}
\int^{t^{n+1}}_{t^n}\hspace{-3mm}
\int^{x_{f}}_{x_{f-1}}\hspace{-3mm}\V s(x,t)\,\d x\,\d t =
\omega_{f-1} \V s_{f-1}^{n+1/2} - \omega_f \V s_f^{\bar t}.
\label{eq:4.6.6}
\end{equation}
This completes the definition of the quadrature rules for the three
integrals on the RHS of (\ref{wf-0}). It remains to extrapolate the solution
$\V w$ to the new quadrature point $(x_f,\bar t)$ near the front.  
%
%
% Section 4.7
%
\subsection{Determination of the physical quantities at $\bar t$}
\label{sect:4.7}
In cases III and IV the margin point passes the cell boundary $x_f$ at
$t^{\ast}$ and goes into the neighbouring cell. The outflow in case
III and the inflow in case IV through the cell boundary at $x_f$ as
well as the source term in the new and old margin cells are
essential for determining the cell average of the margin cells in the
front-tracking method.

In case III the physical quantities flow through the boundary $x_f$ into
the $(f+\halb)^{th}$ cell during the time interval $[t^\ast,\,t^{n+1}]$.
The outflow is approximated by the value at $(x_f,\,\overline t)$,
where $\overline t=\frac{1}{2}(t^{n+1}+t^\ast)$, which is determined
by using Taylor series expansion at the point $(x^n_{Ft},t_n)$
with respect to space and time. Using the
margin-cell-reconstructed slopes and the mass
balance equation (\ref{2.1}), the avalanche depth $h$ at
$(x_f,\,\overline t)$ is then given by
\begin{equation}
\begin{array}{l}
h^{\bar t}_{f}=h_{Ft}^n+\left.\diffx{h}{x}\right|_{(x_{Ft}^{n},t_n)}(x_f-x_{Ft}^{n})
 +\left.\diffx{h}{t}\right|_{(x_{Ft}^{n},t_n)}(\bar t-t^n)
\\*[16pt]\DS\phantom{h^{\bar t}_{f}}
=\left.\diffx{h}{x}\right|_{(x_{Ft}^{n},t_n)}(x_f-x_{Ft}^{n})
 -\left.\diffx{m}{x}\right|_{(x_{Ft}^{n},t_n)}(\bar t-t^n)
\\*[16pt]\DS\phantom{h^{\bar t}_{f}}
=(\sigma^h)^n_f\left(x_f-x^{n}_{Ft}\right)-(\sigma^m)^n_f\left(\bar t-t^n\right),
\\*[16pt]\DS\phantom{h^{\bar t}_{f}}
=(\sigma^h)^n_f\left[\left(x_f-x^{n}_{Ft}\right)-u_{Ft}^n(\bar t-t^n)\right],
\end{array}
\label{hbar}
\end{equation}
where 
$(\sigma^h)^n_f$ and $(\sigma^m)^n_f$ are the slopes of the margin
cell reconstructions defined in \S\,4.1, and we have used the relation
$(\sigma^m)^n_f = u_{Ft}^n (\sigma^h)^n_f$.  Similarly, the momentum
$m^{\bar t}_{f}$ is approximated by
\begin{equation}
m^{\bar t}_{f}
=m_{Ft}^n+\left.\diffx{m}{x}\right|_{(x_{Ft}^{n},t_n)}(x_f-x_{Ft}^{n})
 +\left.\diffx{m}{t}\right|_{(x_{Ft}^{n},t_n)} (\bar t-t^n).
\label{mbar}
\end{equation}
Using the cell reconstruction of the margin $f^{th}$ cell and the
momentum balance equation (\ref{2.2}), equation (\ref{mbar})
becomes
\begin{equation}
\hspace{-1mm}
\begin{array}{l}
m^{\bar t}_{f}
=\left.\diffx{m}{x}\right|_{(x_{Ft}^{n},t_n)}\hspace{-2mm}(x_f-x_{Ft}^{n})
-\left.\diffx{\left(\frac{m^2}{h}
       +\TST\frac{h^2}{2}\beta_x\right)}{x}\right|_{(x_{Ft}^{n},t_n)}
\hspace{-3mm}(\bar t-t^n)
\\*[14pt]\phantom{m^{\bar t}_{f}=m^{\ast}}
+h_{Ft}^n\,s_x(x_{Ft}^{n})\,(\bar t-t^n)
\\*[6pt]\phantom{h^{\bar t}_{f}}
=(\sigma^m)^n_f\bigl(x_f-x^{n}_{Ft}\bigr)
-\left.\biggl(\dfrac{2m}{h}\diffx{m}{x}
-\dfrac{m^2}{h^2}\diffx{h}{x}+\beta_x h \diffx{h}{x}
+\dfrac{h^2}{2}\diffx{\beta_x}{x}\biggr)
\right|_{(x_{Ft}^{n},t_n)}\hspace{-3mm}(\bar t-t^n)
\\*[14pt]
\phantom{h^{\bar t}_{f}}
=(\sigma^m)^n_f\bigl(x_f-x^{n}_{Ft}\bigr)
-\left(2u_{Ft}^n(\sigma^m)^n_f
-\left(u^n_{Ft}\right)^2(\sigma^h)^n_f\right)(\bar t-t^n)
\\*[14pt]
\phantom{h^{\bar t}_{f}}
=u^n_{Ft} \; h^{\bar t}_{f}.
\end{array}
\label{mbar2}
\end{equation}
In case IV the physical quantities at the boundary $(x_f,\,\overline
t)$ are determined in the same way, but the time points are
defined differently: $\overline t=(t^{n+1}+t^\ast)/2$ for case
III and $\overline t =(t^{n}+t^\ast)/2$ for case IV.
%
% Section 4.8
%
\subsection{Summary of the front-tracking algorithm}
\label{sect:4.8}
The front-tracking algorithm may be summarized as follows:
\begin{eqnarray}
\meanW_{f-1/2}^{n+1} &=& \frac12 \meanW_{f-3/4}^{n} + (1 -\alpha_f) \meanW_f^n
- \frac{\overline{\triangle t}}{\triangle x} \V f_f^{\bar t}
+ \frac{\triangle t}{\triangle x} \V f_{f-1}^{n+1/2}
\\
&& + \frac{\omega_{f-1}}{\triangle x} \V s_{f-1}^{n+1/2}
- \frac{\omega_{f}}{\triangle x} \V s_{f}^{\bar t}
\label{eq:4.8.1}
\\[1ex]
%\end{equation}
%\begin{equation}
\meanW_{f+1/2}^{n+1} &=& \alpha_f \meanW_f^n
+ \frac{\overline{\triangle t}}{\triangle x} \V f_f^{\bar t}
+ \frac{\omega_{f}}{\triangle x} \V s_{f}^{\bar t}.
\label{eq:4.8.2}
\end{eqnarray}
Here
\begin{eqnarray}
\alpha_f &=& 
\left\{
\begin{array}{cl}
0 & \textrm{in cases I and III} \\[1ex]
\left(\frac{x_{Ft}^n-x_f}{x_{Ft}^n-x_{f-1/2}}\right)^2 & \textrm{in cases II and IV}
\end{array}
\right.
\label{eq:4.8.3}
\\
\overline{\triangle t} &=&
\left\{
\begin{array}{cl}
0 & \textrm{in case I} \\[1ex]
\triangle t & \textrm{in case II} \\[1ex]
t_{n+1} - t^* & \textrm{in case III} \\[1ex]
t^* - t_n  & \textrm{in case IV}
\end{array}
\right.
\label{eq:4.8.4}
\\
\bar t &=&
\left\{
\begin{array}{cl}
t_n & \textrm{in case I} \\[1ex]
t_{n+1/2} & \textrm{in case II} \\[1ex]
(t_{n+1} + t^*)/2 & \textrm{in case III} \\[1ex]
(t_{n} + t^*)/2  & \textrm{in case IV}
\end{array}
\right.
\label{eq:4.8.5}
\end{eqnarray}
The wheights \( \omega_{f-1}\) and \( \omega_f\) are defined in
(\ref{eq:4.6.3}) and (\ref{eq:4.6.5}). The values of 
\(\V w(x_f,\bar t)\), needed to determine
\(\V f_f^{\bar t}\) and \(\V s_f^{\bar t}\), are defined 
in (\ref{hbar}) and (\ref{mbar2}). 
This completes the definition of the update at the front margin.
The tail can be treated completely analogously.
%
% Section 5
%
\section{Numerical Results}
%
%
% Section 5.1
%
\subsection{Travelling shock wave}
In this test problem we are concerned with granular flow on a plane
($\lambda \kappa = 0$) inclined chute ($0\le x \le 36$
dimensionless units), where the
internal and basal friction angles are both presumed to be equal to
the inclination angle, $\phi=\delta=\zeta=40^\circ$. That implies a
non-accelerative flow, $s_x=0$, whose earth pressure coefficient is
constant $K_x=K_{x_{act}}=K_{x_{pass}}$. Selecting $\varepsilon=1$ and
using (\ref{2.5}) yields $\beta_x=\varepsilon\cosZ K_x=1.84477$. A
jump of thickness $H=h^-/h^+=3$ with $h^+=0.3,\ h^-=0.9$ is presumed
at $x=24$.  By virtue of (\ref{VelDiff2}) the velocity difference is
then determined, $u^+-u^-=1.2148317$, where the positive sign is
selected.  Since an instability was expected close to $u=0$ as a
singularity by sgn$(u)$, the downslope velocity is assumed to be
$u^-=0.1$, so that the term sgn$(u)$ is always unity.  The initial
condition of this test problem is defined as follows:
\begin{align}
%\begin{array}{l}
h(x,\,0)&=\left\{
\begin{array}{l}
0.3, \quad \hbox{for}\quad 0\leq x <24,\\
0.9, \quad \hbox{for}\quad 24\leq x<36,
\end{array}
\right. \\[12pt]
u(x,\,0)&=\left\{
\begin{array}{l}
1.3148317, \quad \hbox{for}\quad 0\leq x <24,\\
0.1, \phantom{148317}\quad \hbox{for}\quad 24\leq x\leq 36.
\end{array}
\right. 
%\end{array}
\label{8.2}
\end{align}
From (\ref{vel_jump}) the velocity of the upslope travelling wave is
then expected as $V_n=-0.50741585$. For the boundary condition a
constant inflow at $x=0$ and an outflow condition for $x=36$ are
introduced.

%
% Section 5.1.1
%
\subsubsection{Lagrangian technique}

By the Lagrangian moving grid method the governing equations
(\ref{2.1}) and (\ref{2.2_1}) are solved by virtue of
(\ref{3.2})--(\ref{3.5}). The initial depth, $h_j^0$, of the $j^{th}$
element is taken to be the cell average of the exact initial profile.
The initial velocity of the boundary, $u_j^0$, is given by the volume
weighted velocity of the adjacent cells. They are
\begin{equation}
\DS h_j^0=\dfrac{\DS\int^{b^0_{j}}_{b^0_{j-1}}\hspace{-2mm}h(x,\,0)\,\d x}%
{\DS{b_{j}^0}-{b_{j-1}^0}},\qquad
\DS u_j^0=\dfrac{\DS\int^{c^0_{j+1}}_{c^0_j}\hspace{-3mm}h(x,\,0)\,u(x,\,0)\,\d x}%
{\DS\int^{c^0_{j+1}}_{c^0_j}\hspace{-4mm}h(x,\,0)\,\d x},
\label{Lag_ini}
\end{equation}
where $b_{j}^0$ and ${b_{j-1}^0}$ are the boundaries of the $j^{th}$
cell at $t=0$, and $c^0_j$ denotes the initial centre of the $j^{th}$
cell.
\begin{figure}[t!]
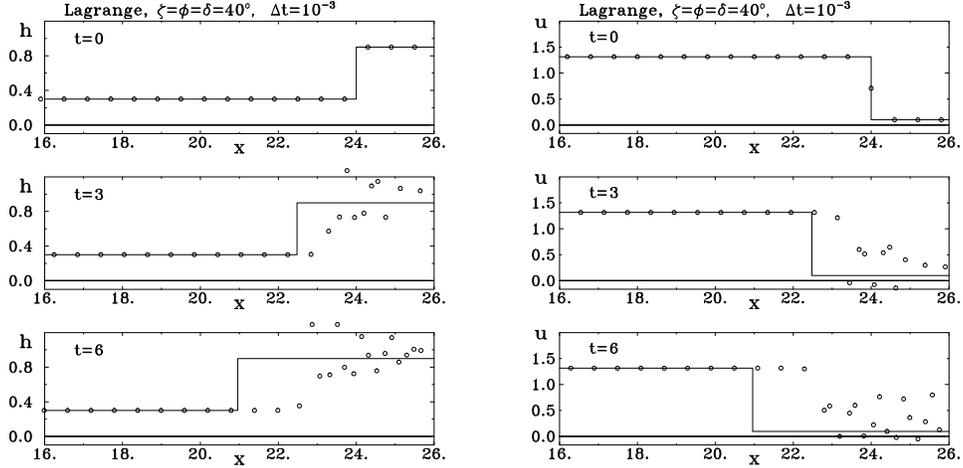

%\vskip4mm
  \begin{center}
    \includegraphics[angle=0.0,scale=0.48]{figures/work/gibbs/tai/Diss2/image/LaH.epsi}\hfill
    \includegraphics[angle=0.0,scale=0.48]{figures/work/gibbs/tai/Diss2/image/LaU.epsi}
\vskip12pt
  \caption[{Figjump}]{Depth (left) and the corresponding 
    velocity (right) profiles of the upslope travelling wave at
    $t=0,\,3,\,6$, where circles denote the computed results at the
    cell centres and the solid line indicates the exact solution. The
    time step is taken to be $\dt=10^{-3}$ dimensionless time unit.
    \label{LJump}}
  \end{center}
\end{figure}

The constant inflow and outflow boundary conditions are executed by
setting the depth gradient $\partial h/\partial x$ at $x = b_0(t)$ and
$x = b_N(t)$ equal to zero, so that $\d u/\d
t=0\:\Rightarrow\:u_0(t)=u_0(0)$ and $u_N(t)=u_N(0)$ for $t>0$ because
the flow is on an non-accelerating slope $s_x=0$.

Fig.\,\ref{LJump} demonstrates the simulated results ($N=60$);
oscillations develop as the shock wave passes through, and these
persist even if the time step is selected to be very small. The
velocities of the cell boundary after the shock are sometimes faster
or slower than they should be and therefore oscillations take place.
These oscillations propagate downslope as time increases and no shock
wave propagates upslope. This indicates that the Lagrangian moving
grid technique is ill behaved and cannot describe the travelling shock
wave.
%
% Section 5.1.2
%
\subsubsection{Eulerian shock-capturing methods}
%
%
% 6.9
%
\begin{figure}[t!]
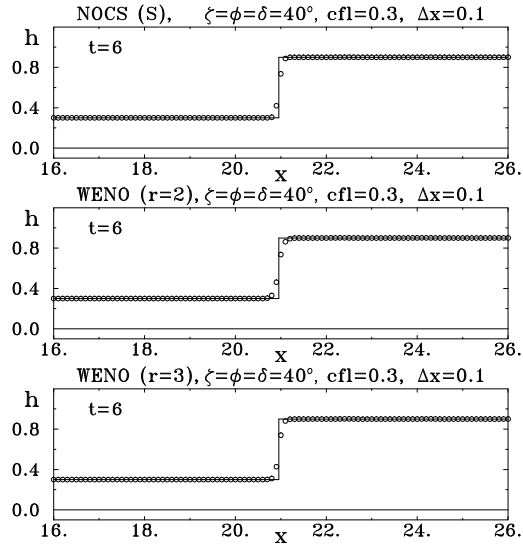

  \begin{center}
    \includegraphics[angle=0.0,scale=0.56]{figures/work/gibbs/tai/Diss2/image/einNSH.epsi}
    \includegraphics[angle=0.0,scale=0.56]{figures/work/gibbs/tai/Diss2/image/einWr2H.epsi}
    \includegraphics[angle=0.0,scale=0.56]{figures/work/gibbs/tai/Diss2/image/einWr3H.epsi}
\vskip12pt
  \caption[{Figjump}]{Depth  profiles of the upslope 
    travelling wave computed by the NOC scheme at $t=6$ with $N=360$.
    The solid lines indicate the exact solution and circles mean the
    computed results.
    \label{WcomJump}}
  \end{center}
\end{figure}

The NOC scheme is applied to (\ref{2.1}) and (\ref{2.2}) on a 1D grid
with $90$ and $360$ gridpoints, respectively.  The initial conditions
are transferred to the mean values over the cells before the computing
commences,
\begin{equation}
\DS h^0_j=\dfrac{1}{\dx}\int^{x_{j+1/2}}_{x_{j-1/2}}h(x,\,0)\d x,
\qquad
\DS u^0_j=\dfrac{1}{h^0_j}\int^{x_{j+1/2}}_{x_{j-1/2}}h(x,\,0)u(x,\,0)\d x.
\label{cellaverge}
\end{equation}
The constant inflow boundary condition is implemented by the
assignments $h_0(t)=h_0(0)$ and $m_0(t)=m_0(0)$ at $x=0$. The outflow
boundary condition is described by setting $\dd hx=0$ and $\dd mx=0$
at $x=36$, where they are
\begin{equation}
U_N=\left(4\,U_{N-1}-U_{N-2}\right)/3,\quad \hbox{for}\quad U=h,~m,
\end{equation}
by using the cell averages of the closest cells for a second
order extrapolation.

Three different cell reconstructions were tested: the NOC scheme with
superbee limiter (NOCS-S), piecewise linear (r=2) and quadratic
(r=3) WENO reconstructions \cite{LevyPuppoRusso:1998}. Fig.\,\ref{WcomJump}
demonstrates the simulated avalanche depth of the
travelling wave problem (circles) and a comparison with the exact
solution (solid line) at $t=6$ dimensionless time units. All of them
are able to adequately describe this travelling shock wave problem.
%
% Section 5.2
%
\subsection{Parabolic similarity solution}
This section is concerned with the simulation of the parabolic
similarity solution outlined in \S\,2.2.  In the test problem the
parabolic avalanche body is considered to slide on an inclined flat
plane in the domain $0\le x \le 36$ dimensionless length units with
constant inclination angle $\zeta=40^\circ$. The basal and internal
friction angles are simultaneously selected to be $30^\circ$, and the
initial condition is chosen to be $g_0=1$ and $p_0=0$. On the inclined
plane the initial depth and velocity distributions are mapped into
\begin{equation}
\left\{
\begin{array}{l}
h(x,\,0)=1-\bigl((x-4)/3.2\bigr)^2,\\[10pt]
u(x,\,0)=1.2,
\end{array}
\right.
\quad\hbox{for}~x\in[0.8,\,7.2].
\end{equation}
Our choice of the initial velocity, $u(x,\,0)=u_0=1.2$, guarantees that
condition (\ref{Vcon}) will be satisfied for all times.  This problem will
serve as the standard test problem for the resolution of the depth
profile and the determination of the margin locations.

%
% Section 5.2.1
%
\subsubsection{Lagrangian technique}

In the Lagrangian moving grid technique the model equations
(\ref{2.1}) and (\ref{2.2_1}) are solved by virtue of
(\ref{3.2})--(\ref{3.5}) on a 1D grid. The boundary condition is given
by setting the heights at the margin (front and rear) points to be
equal to zero, $h_0(x,\,t)=0$ and $h_N(x,\,t)=0$.

\begin{figure}[t!]
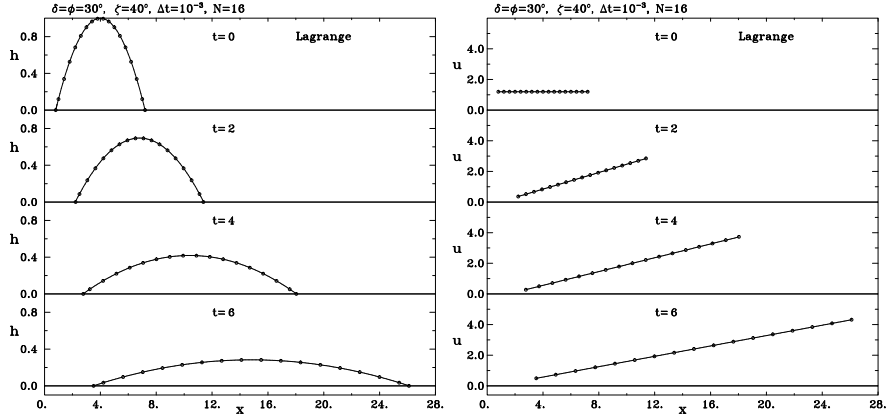

%\vskip4mm
  \begin{center}
    \includegraphics[angle=0.0,scale=0.34]{figures/work/gibbs/tai/Lawine/1D/Lagrange/Sim-data/lag16H.epsi}
    \includegraphics[angle=0.0,scale=0.34]{figures/work/gibbs/tai/Lawine/1D/Lagrange/Sim-data/lag16U.epsi}
\vskip12pt
  \caption[{Figjump}]{Depth (left) and the corresponding 
    velocity (right) profiles of the parabolic similarity solution
    (Problem I) computed by the Lagrangian moving-grid scheme at the
    dimensionless time units $t=0,\,2,\,4,\,6$, where the avalanche
    body is divided into $16$ cells, and the time interval is
    $\dt=10^{-3}$.\label{Lag16}}
  \end{center}
\end{figure}

Fig.\,\ref{Lag16} illustrates the simulated result at the
dimensionless time units $t=0,\,2,\,4,\,6$ with cell number $N=16$, in
which the circles denote the computed results and the solid line
indicates the exact solution. The avalanche body extends as it flows
down and still keeps the parabolic depth profile. The velocity is
keeping a linear distribution through the bulk body. It ensures the
symmetric depth profile during the motion.

\begin{figure}[t!]
%\vskip4mm
  \begin{center}
    \includegraphics[angle=0.0,scale=0.7]{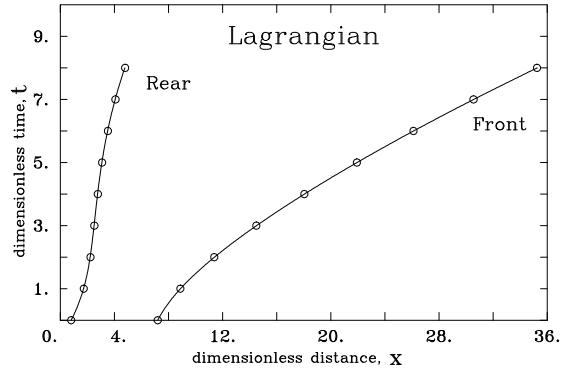}
\vskip12pt
  \caption[{Figjump}]{Locations of the front and
    rear edges of the avalanche body in the parabolic similarity
    solution problem as they evolve in time. The circles denote the
    computed results by the Lagrangian moving grid technique ($N=16$),
    and the solid lines indicate the exact margin positions. They are
    in excellent agreement.\label{Lagrand}}
  \end{center}
\end{figure}

From the simulated results it follows that the Lagrangian moving grid
technique can not only describe the depth profile well but also
determines the margin locations of the similarity solution very
accurately. There is excellent agreement between the simulated results
and the exact solutions, see Fig.\,\ref{Lag16}.  The motions of the
front and rear edges of the avalanche body in the similarity solution
are illustrated in Fig.\,\ref{Lagrand}. The circles denote the
computed results by the Lagrangian moving grid technique and the solid
lines indicate the exact locations of the margins. They are also in
excellent agreement.

% Fig. 8.1
\begin{figure}[t!]
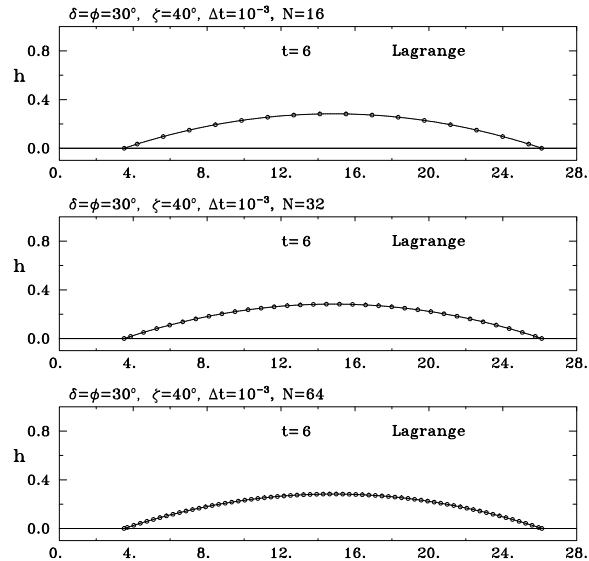

%\vskip4mm
  \begin{center}
    \includegraphics[angle=0.0,scale=0.45]{figures/work/gibbs/tai/Lawine/1D/Lagrange/Sim-data/La16e.epsi}\\[5pt]
    \includegraphics[angle=0.0,scale=0.45]{figures/work/gibbs/tai/Lawine/1D/Lagrange/Sim-data/La32e.epsi}\\[5pt]
    \includegraphics[angle=0.0,scale=0.45]{figures/work/gibbs/tai/Lawine/1D/Lagrange/Sim-data/La64e.epsi}
\vskip12pt
  \caption[{Figjump}]{Depth profiles
    computed by the Lagrangian moving-grid technique for the parabolic
    similarity solution problem (Problem I), where the avalanche body
    is divided into different numbers of cells, $N=16,\,32,\,64$.  All
    the results are shown at $t=6$ and the computational time interval
    is $\dt=10^{-3}$. The number of the cells does not influence the
    good agreement between the simulated results (circles) and the
    exact solutions (solid line).\label{Lagcom}}
  \end{center}
\end{figure}

The Lagrangian method is also tested by different grid numbers.
Fig.\,\ref{Lagcom} shows the results computed by different grid
numbers, $N=16,\,32$, and $64$, respectively. With different grid
numbers this method can always keep the excellent resolutions when
compared with the exact solutions.

Caculations were also performed with initial condition $p_0\neq 0$;
results turned out to be similarly convincing as the above ones. For
this reason they are not presented here \cite{Tai:2000}.
% Fig. 8.2
%
%Fig. 8.3
%
%\clearpage
%
% Section 5.2.2
%
\subsubsection{Eulerian technique}

In \S\,3.2, the Eulerian schemes are based on the model equations
(\ref{2.1}) and (\ref{2.2}) in conservative form, so that the velocity
outside the avalanche body (inclusive the margin point) is not
defined. Intuitively, adding a thin layer of the material over the
whole computational domain could be used to treat the grain free
regions. Another trick can also be introduced, in which all the physical
variables are set to zero if $h=0$. This would be reasonable since
$h=0~\rightarrow~m=hu=0$.

\begin{figure}[t!]
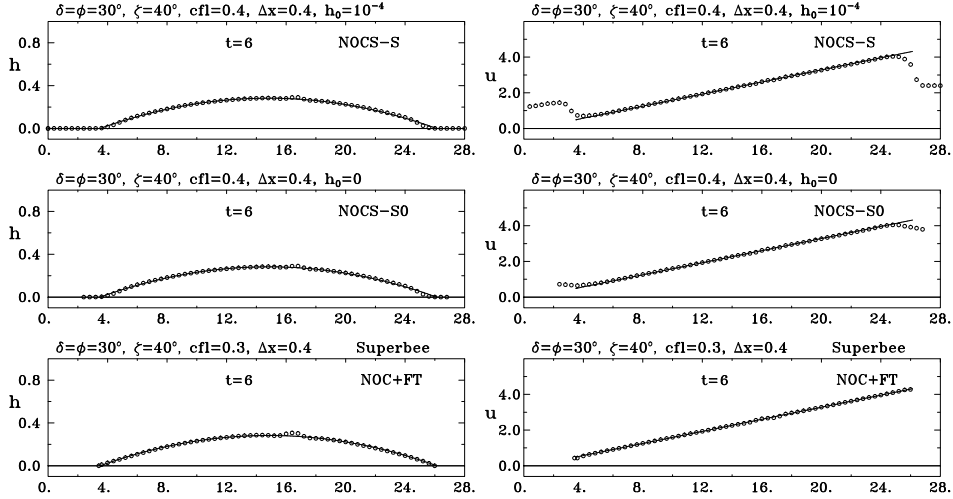

%\vskip4mm
  \begin{center}
    \includegraphics[angle=0.0,scale=0.44]{figures/work/gibbs/tai/Lawine/1D/NOCS/Sim-data/nocs90sHe.epsi}
    \includegraphics[angle=0.0,scale=0.44]{figures/work/gibbs/tai/Lawine/1D/NOCS/Sim-data/nocs90sUe.epsi}\\[5pt]
    \includegraphics[angle=0.0,scale=0.44]{figures/work/gibbs/tai/Lawine/1D/NOCS/Sim-data/nocs900He.epsi}
    \includegraphics[angle=0.0,scale=0.44]{figures/work/gibbs/tai/Lawine/1D/NOCS/Sim-data/nocs900Ue.epsi}\\[5pt]
    \includegraphics[angle=0.0,scale=0.44]{figures/work/gibbs/tai/Lawine/1D/NOCS/Sim-data/noel90sHe.epsi}
    \includegraphics[angle=0.0,scale=0.44]{figures/work/gibbs/tai/Lawine/1D/NOCS/Sim-data/noel90sUe.epsi}
\vskip12pt
  \caption[{Figjump}]{
    Depth (left) and velocity (right) profiles of the parabolic
    similarity solution computed by the NOC scheme with Superbee
    limiter. In the top panels, a thin layer with $h_0=10^{-4}$ is
    added to the whole computational domain.  In the middle panels,
    all physical variables are set to zero if $h=0$. Whilst, the
    bottom pannels demonstrate the results from the scheme with
    front-tracking method. The whole computational domain is divided
    into $90$ cells ($N=90$), the circles denote simulated results
    and the solid lines represent the exact solution. The results
    show that the added thin layer does not influence the depth
    profile very much, if it is sufficiently small, but the margin
    locations can not be exactly determined without the front-tracking
    method. An oscillation near the middle of the avalanche
    (local maximum) is visible in all three calculations.
    \label{NOCSthin}}
  \end{center}
\end{figure}

Fig.\,\ref{NOCSthin} illustrates the
comparison between the computed results obtained from the NOC scheme,
where a thin layer $h_0=10^{-4}$ repectively $h_0=0$ is added
over the whole computational domain, and from the scheme with our
front-tracking method.
All the three results of the depth profiles are acceptable except for
the oscillation near the top. However, have a look at the velocity
profiles in these figures, there are several cells with $\dd ux<0$
around the margins. This violates the assumption $\dd ux >0$ in the
parabolic similarity solution problem.
Moreover, the results show that there is large numerical
diffusion around the margins (i.e. the margins move further than they
should do) without the front-tracking method.
For both reasons, the front-tracking method is needed
to determine the location of the margins.

% Fig. 8.6
%

%By inspecting the results in Fig.\,\ref{NOCSthin} it is found that the
%thickness of the added thin layer (either $h_0=10^{-4}$ or $h_0=0$)
%does not have an obvious influence on the results of the depth
%profile.  But there are large velocity gradients around the margins.
%It follows from the momentum balance equation (\ref{2.2}), since the
%regions outside of the margin are covered by the uniformly added thin
%layer, that there is no acceleration due to the depth gradient, $\dd hx$, in
%the momentum balance equation. So, the velocity in this region will
%tend to have the value of the moving co-ordinate, $u_0(t)$, as defined
%in (\ref{2.17}). On the other hand, inside the avalanche body there is
%a permanent contribution from the depth gradient in the momentum
%balance equation. Therefore a jump of velocity develops around the margin.
%Further, the wrong cell reconstruction over the margin cells, see \S\,4.1,
%causes the unphysical diffusion, which results in the dislocation
%of the margin.
% Fig. 8.7
\begin{figure}[t!]
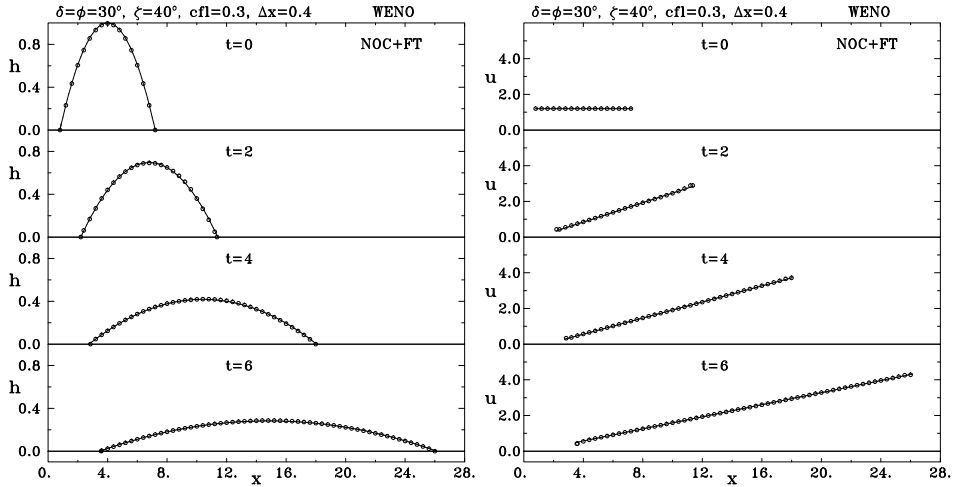

%\vskip4mm
  \begin{center}
    \includegraphics[angle=0.0,scale=0.44]{figures/work/gibbs/tai/Lawine/1D/NOCS/Sim-data/noel90wH.epsi}    
    \includegraphics[angle=0.0,scale=0.44]{figures/work/gibbs/tai/Lawine/1D/NOCS/Sim-data/noel90wU.epsi}    
\vskip12pt
  \caption[{Figjump}]{
    Depth (left) and velocity (right) profiles of the parabolic
    similarity solution at $t=6$ computed by the NOC scheme with
    front-tracking and piecewise quadratic WENO cell reconstruction.
    The whole computational domain is divided into
    $90$ cells ($N=90$) and the Courant number is selected to be
    $0.3$.  The margin locations are well described and the oscillation
    near the center is successfully removed.
    \label{NFTcom}}
\end{center}
\end{figure}
%
%\vspace{8mm}
% Fig. 8.11
\begin{figure}[h!]
%\vskip4mm
  \begin{center}
    \includegraphics[angle=0.0,scale=0.7]{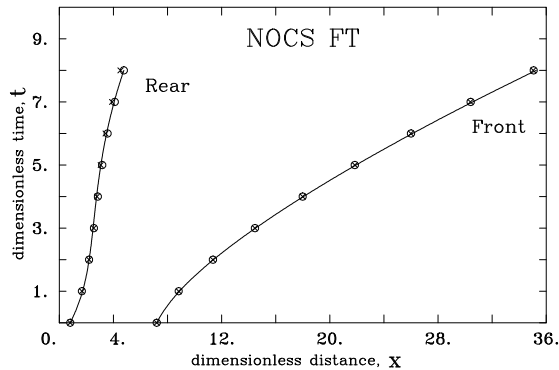}
\vskip12pt
  \caption[{Figjump}]{Front and rear edges of the avalanche body in
    the parabolic similarity solution simulated by the
    NOC-front-tracking scheme as they evolve in time.
    \lq\lq$\circ$\rq\rq\ denotes the computed results obtained with
    the piecewise quadratic WENO cell reconstruction, \lq\lq$\times$\rq\rq\ 
    means the results deduced with Superbee limiter and solid lines
    indicate the exact margin solution.\label{NOCrand}}
  \end{center}
\end{figure}

Let us discuss the origin of the oscillation near the center of the avalanche.
When one recomputes the solution using unlimited central differences
for $\meanW^{\,\prime}$, the oscillation disappears.
Therefore, we have the following paradoxical
situation: the introduction of TVD-limiters, which are needed to stabilize
the solution in the presence of discontinuities, may destabilize the
solutions in smooth regions! In fact, this is not entirely surprising,
since in the presence of limiters the fluxes depend only Lipschitz-continuously
on the data.

We have therefore experimented with more smooth reconstructions,
namely piecewise quadratic WENO interpolants of Jiang and Shu 
\cite{JiangShu:1996} and Levy, Puppo and Russo 
\cite{LevyPuppoRusso:1998}, which depend smoothly on the data
and are at the same time nonoscillatory at discontinuities.
In the margin cells, we still use the piecewise linear reconstructions
introduced in Section \ref{sect:4.1}, and in the two cells adjacent to
the margin cells, we use a piecewise linear WENO reconstruction.
We have experimented both with second- and third-order quadrature
rules in time. In our experience, both yield comparable results.
Fig.\,\ref{NFTcom} demonstrates the results for these reconstructions
combined with our front-tracking method. The margin
locations are well described by the front-tracking method,
and the oscillation near the center is successfully removed
(compare the bottom pannels in Fig.\,\ref{NOCSthin} and
Fig.\,\ref{NFTcom}).

Fig.\,\ref{NOCrand} shows the computed front and rear edges of the
avalanche body in the parabolic similarity solution 
as they evolve in time.  \lq\lq$\circ$\rq\rq\ denotes the
computed results obtained by the NOC scheme with the piecewise quadratic WENO cell
reconstruction, \lq\lq$\times$\rq\rq\ means the results deduced with
Superbee limiter and solid lines indicate the exact margin solution.
Both the Superbee limiter and the piecewise quadratic WENO cell reconstruction
for the NOC front-tracking schemes can yield good agreement of the
determined margin locations with the exact solutions.
\begin{table}[t!]
\label{tab8_1}
\small
\begin{center}
\caption{Error (\ref{errorMeasure2}) of the different schemes.}
\begin{tabular}{c|cccc|cc}
\\[-8pt]
 $$ & NFT(S90) & NFT(W90) 
 & NFT(S360) & NFT(W360) & Lag(16) & Lag(32)
\rule{0in}{2ex}\\[.5ex] \hline\\[-8pt]
   &  &   \multicolumn{2}{c}{$E~(\times 10^{-3})$} 
   &  &   \multicolumn{2}{c}{$E~(\times 10^{-3})$} 
\rule{0in}{2ex}\\[.5ex] \hline\\[-8pt]
  $t=1$ & $ 9.0247$ & $7.6096$ & $0.8816$ & $0.8813$ & $1.7130$ & $0.2937$ 
\rule{0in}{2ex}\\[.5ex] 
  $t=2$ & $11.3333$ & $7.5679$ & $0.8532$ & $0.9024$ & $1.7764$ & $0.3664$ 
\rule{0in}{2ex}\\[.5ex] 
  $t=3$ & $12.8854$ & $7.1051$ & $0.7336$ & $0.9492$ & $1.8944$ & $0.4135$ 
\rule{0in}{2ex}\\[.5ex] 
  $t=4$ & $16.3450$ & $6.0860$ & $0.8484$ & $1.0875$ & $1.8888$ & $0.4413$ 
\rule{0in}{2ex}\\[.5ex] 
  $t=5$ & $17.4274$ & $6.5434$ & $0.8672$ & $1.1203$ & $1.8974$ & $0.4492$ 
\rule{0in}{2ex}\\[.5ex] 
  $t=6$ & $18.8426$ & $6.5340$ & $1.1109$ & $1.3203$ & $1.9474$ & $0.4658$ 
\rule{0in}{2ex}\\[.5ex] 
  $t=7$ & $16.5043$ & $6.0608$ & $1.3970$ & $1.0435$ & $1.8817$ & $0.5026$ 
\rule{0in}{2ex}\\[.5ex] 
  $t=8$ & $15.3966$ & $6.0762$ & $1.8558$ & $1.1981$ & $1.9526$ & $0.4830$ 
\rule{0in}{2ex}\\[.5ex] 
\hline 
\end{tabular}
%\begin{minipage}{19cm}
%\end{minipage}
\end{center}
\end{table}
%\end{landscape}

The use of the Superbee limiter results in a small delay of the
avalanche body, i.e. a slower velocity both at the front and the rear.
The reason is that the Superbee limiter tends to be overcompressive
in smooth regions of the solution, and therefore it does not give the appropriate
flux at the boundaries between the internal and the margin cells.

In order to obtain some quantitative information on the accuracy
of the schemes, we introduce an error measure for the depth by
\begin{equation}
E=\dfrac{\sum\limits^N_{j=0} \left|h_j-\overline h^{exact}_j\right|}%
{\sum\limits^N_{j=0} \overline h^{exact}_j};
\label{errorMeasure2} 
\end{equation}
where $\overline
h^{exact}_j$ denotes the $j^{th}$ cell averaged depth of the exact
solution. The errors of the Lagrangian method, the NOC front-tracking
scheme with Superbee limiter (NFT(S)) and piecewise quadratic WENO
interpolations (NFT(W)) at $t=1$ to $8$ dimensionless time units are
shown in Table.  Here, the Eulerian schemes are tested by using $N=90$
and $N=360$, respectively, and for the Lagrangian scheme $N=16$ and
$N=32$ are used. The Lagrangian method results in the least errors,
obviously smaller than all the Eulerian schemes. It also converges
at a better rate.
%
% Section 5.3
%
\subsection{Upward Moving Shock Wave}
Shock formations are often observed when the avalanche slides into the
run-out horizontal zone. Here the front part comes to rest, while the
tail accelerates further and its velocity becomes supercritical.
In \cite{TaiNoelleGrayHutter:2000} a comparison was made between our shock-capturing
method and the Lagrangian moving grid technique for the case
of coinciding basal and internal friction angles.
Here we compute a flow with basal friction angle $\phi = 38^\circ$
and internal friction angle $\delta = 35^\circ$. As a consequence,
we have a jump in the earth-pressure coefficient $K_x$ when the
flow changes from an expanding ($u_x > 0$) to a contracting region
($u_x < 0$).
%
% Figure Exp.3
%
\begin{figure}[h!]
\begin{center}
\includegraphics[angle=0.0,scale=0.5]{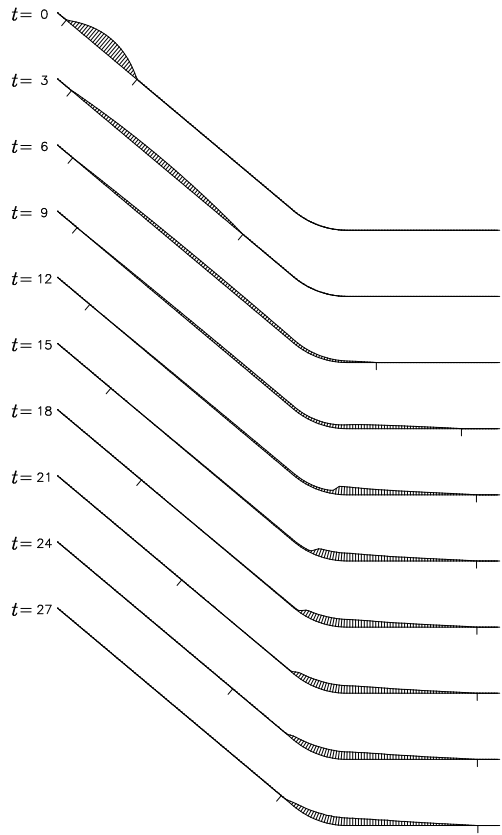}
\vskip12pt
    \caption[{NOCS_Exp3}]{
      Process of the avalanche simulated by the shock-capturing
      and front-tracking NOC
      method at $t=0,\,3,\,6,\dots,\,27$ dimensionless time units.
      As the front reaches the run-out zone and comes to rest, the
      part of the tail accelerates further and the avalanche body
      contracts.  Once the velocity becomes supercritical, a shock
      wave develops, which moves upward.
      The dashes below the graphs mark the tail and the head
      of the avalanche. \label{NOCS_Exp3}}
\end{center}
\end{figure}

The setup is as follows: 
The granular material released from a parabolic cap slides down an
inclined plane and merges into the run-out horizontal zone. The centre
of the cap is initially located at $x=4.0$ and the initial radius and
the height are $3.2$ and $1.0$ dimensionless length units,
respectively. The inclination angle of the inclined plane is
$40^\circ$ and the (linear and continuous) transition region
lies between $x=21.5$ and $x=25.5$. We use 180 gridpoints
and a CFL-number of 0.4.

Fig.\,\ref{NOCS_Exp3} illustrates the simulated process as the
avalanche slides on the inclined plane into the horizontal run-out
zone (so initially the flow is expanding).
The avalanche body extends on the inclined plane until the
front reaches the run-out zone. Here the basal friction is
enough to bring the front of the granular material to rest while the
rear part accelerates further. Therefore, the flow becomes
contracting in the transition zone. At this stage, a shock (surge)
wave is created ($t=12$), which moves upward. Such shock waves make
the Lagrangian method unstable, if no artificial
viscosity is applied (see \cite{TaiNoelleGrayHutter:2000}).
Our non-oscillatory central front-tracking scheme handles
both the shock wave and the margins of the avalanche well.

%
% Section 6
%
\section{Conclusion}

In this paper we have developed a Lagrangian and an Eulerian
shock-capturing finite-difference scheme with front-tracking
for the spatially one-dimensional Savage-Hutter equations of granular
avalanches. The purpose was to reproduce the temporal evolution of the
avalanche geometry and downslope velocity under situations when
internal shocks may occur. This happens e.g. when an avalanche of
finite mass moves from an inclined chute into the horizontal run-out
zone and, in the transition zone, is deccelerated from a supercritical
flow state to a subcritical state.
The Lagrangian scheme (which is excellent for smooth solutions)
develops unphysical oscillations when the solution contains,
or develops, shock discontinuities.
In order to compute discontinuous solutions, we propose to use a
conservative shock-capturing finite difference scheme.
We adapt the second-order accurate staggered scheme of
Nessyahu and Tadmor \cite{NessyahuTadmor:1990} to the Savage-Hutter equations.
The staggered approach avoids the use of characteristic decompositions which
are needed in standard upwind schemes, but are not known for
the Savage-Hutter equations.
We show that our non-oscillatory central (NOC) scheme
reproduces both smooth and shock solutions adequately
except for the following two problems: First, oscillations may occur near smooth
extrema due to the presence of piecewise linear reconstructions
with TVD-type limiters. These oscillations disappear when one uses
piecewise quadratic cell reconstructions in the interior of the avalanche.
Second, our NOC scheme (and in fact, any Eulerian scheme) does not
capture the vacuum-boundary accurately. This may lead to serious stability problems.
We improve the treatment of the free boundary by combining the scheme
with a front-tracking method applied to the margin
cells. In the spirit of the Nessyahu-Tadmor scheme, we do
not make use of the vacuum Riemann-problem, but rely on
a new piecewise linear reconstruction at the vacuum boundary
and carefully chosen Taylor-extrapolations for the corresponding
numerical fluxes.
With such a combination of an internal Eulerian NOC scheme and
a Lagrangian \lq\lq boundary scheme\rq\rq\ two standard test problems -- an
upward moving shock and a parabolic cap moving down an inclined plane
-- could be well reproduced.
The scheme also produces satisfactory results for the more realistic problem 
mentioned above: an avalanche moving down an inclined plane and coming to rest 
at a flat run-out. Here an upward moving shock wave develops from smooth
data, and the flow changes from being expanding to contracting 
ahead of the shock. In this situation, the earth pressure coefficient 
changes discontinuously, so we are facing the full difficulties inherent
in the Savage-Hutter model.

Several questions remain and await further study:
\begin{itemize}
  
\item The shock-capturing NOC numerical method including the
  front-tracking scheme should be extended to two-dimensional flows.
  This is work under progress.
  
\item The original Lagrangian moving grid scheme could also be
  developed as a shock-capturing scheme. Here the main difficulty
  would be in the determination of the correct grid velocity.

\end{itemize}
We are working on these topics and will report on results in due time.

\begin{acknowledgment}
  Y.C. Tai, J.M.N.T. Gray and K. Hutter acknowledge financial support
  from the Deutsche Forschungsgemeinschaft via SFB 298, \lq\lq
  Deformation und Versagen metallischer und granularer
  Strukturen\rq\rq at Darmstadt University of Technology.
  S. Noelle was supported by SFB 256, \lq\lq
  Nichtlineare Partielle Differentialgleichungen\rq\rq, 
  at Bonn University.
\end{acknowledgment}
%%% Local Variables: 
%%% mode: latex
%%% TeX-master: t
%%% End: 

%% End of article:

%% This command is necessary! ==>>
\end{article}
\end{document}